\documentclass[12pt, a4paper]{amsart}

\usepackage{geometry}
\usepackage{amsmath, amssymb}
\usepackage{enumerate}
\usepackage{graphicx}
\graphicspath{{fig/}}

\newtheorem{remark}{Remark}
\newtheorem{theorem}{Theorem}
\newtheorem{algorithm}{Algorithm}

\makeatletter %
\@addtoreset{equation}{section} %
\makeatother

\DeclareMathOperator*{\argmin}{argmin}

\newcommand{\mca}{\mathcal{A}}
\newcommand{\mcr}{\mathcal{R}}
\newcommand{\mct}{\mathcal{T}}

\newcommand{\br}{\mathbb{R}}
\newcommand{\rd}{\mathrm{d}}
\newcommand{\veps}{\varepsilon}

\newcommand{\bn}{\boldsymbol{n}}
\newcommand{\boeta}{\boldsymbol{\eta}}
\newcommand{\bx}{\boldsymbol{x}}
\newcommand{\bxi}{\boldsymbol{\xi}}
\newcommand{\bzeta}{\boldsymbol{\zeta}}

\newcommand{\abs}[1]{\left\vert#1\right\vert}
\newcommand{\brac}[1]{\left(#1\right)}
\newcommand{\norm}[1]{\left\Vert#1\right\Vert}

\begin{document}

\title[The auxiliary function method]{The auxiliary function method for waveform based earthquake location}

\author{Hao Wu}
\address{Department of Mathematical Sciences \\
  Tsinghua University \\
  Beijing, China 100084 \\
  email: hwu@tsinghua.edu.cn }
  
\author{Jing Chen}
\address{Department of Mathematical Sciences \\
  Tsinghua University \\
  Beijing, China 100084 \\
  email: jing-che16@mails.tsinghua.edu.cn }

\author{Hao Jing}
\address{Department of Mathematical Sciences \\
  Tsinghua University \\
  Beijing, China 100084 \\
  email: jinghao0320@gmail.com }
  
  \author{Ping Tong}
\address{College of Science \\
  Nanyang Technological University, Singapore \\
  email: tongping@ntu.edu.sg }

\author{Dinghui Yang}
\address{Department of Mathematical Sciences \\
  Tsinghua University \\
  Beijing, China 100084 \\
  email: dhyang@math.tsinghua.edu.cn }

\date{\today}

\begin{abstract}
	This paper introduces the auxiliary function method (AFM), a novel, fast and simple approach for waveform based earthquake location. From any initial hypocenter and origin time, we can construct the auxiliary function, whose zero set contains the real earthquake hypocenter and the origin time. In most of situations, there are very few elements in this set. The overall computational cost of the AFM is significantly less than that of the iterative methods. According to our numerical tests, even for large noise, the method can still achieve good location results. These allow us to determine the earthquake hypocenter and the origin time extremely fast and accurate. \medskip
	
\end{abstract}

\maketitle


\section{Introduction} \label{sec:intro}

In this work, we present a novel approach to solve the inverse problem \cite{WuChHuYa:16} to determine the real earthquake hypocenter $\bxi_T$ and the origin time $\tau_T$
\begin{equation} \label{eqn:optim}
	(\bxi_T,\tau_T)=\argmin_{\bxi,\tau}\sum_{r\in\mcr}\chi_r(\bxi,\tau),
\end{equation}
in which $\chi_r(\bxi,\tau)$ is the misfit function
\begin{equation} \label{fun:mis}
	\chi_r(\bxi,\tau)=\frac{\int_0^T\abs{d_r(t)-s(\boeta_r,t)}^2\rd t}{2\int_0^T\abs{d_r(t)}^2\rd t},
\end{equation}
for the $r-$th receiver. 
The time-series $d_r(t)$ is the real earthquake signal which occurred at $(\bxi_T,\tau_T)$ and was recorded at  receiver $r$. The set $\mcr$ contains all of the receivers that we use for inversion. And the synthetic earthquake signal $s(\bx,t)$ is corresponding to the initial hypocenter $\bxi$ and the initial origin time $\tau$. For model simplicity, they can be regarded as the solutions
\begin{equation} \label{eqn:rel_ds}
	d_r(t)=u(\boeta_r,t;\bxi_T,\tau_T), \quad s(\bx,t)=u(\bx,t;\bxi,\tau),
\end{equation}
of the following acoustic wave equation
\begin{equation} \label{eqn:wave}
	\frac{\partial^2 u(\bx,t;\bxi,\tau)}{\partial t^2}=\nabla\cdot\brac{c^2(\bx)\nabla u(\bx,t;\bxi,\tau)} \\
		+f(t-\tau)\delta(\bx-\bxi), \quad \bx,\bxi\in\Omega,
\end{equation}
with initial-boundary conditions
\begin{align} 
	& u(\bx,0;\bxi,\tau)=\partial_t u(\bx,0;\bxi,\tau)=0, \quad \bx\in\Omega, \label{ic:wave} \\
	& \bn\cdot \brac{c^2(\bx)\nabla u(\bx,t;\bxi,\tau)}=0, \quad \bx\in\partial \Omega. \label{bc:wave}
\end{align}
Here $\boeta_r$ denotes the location of the $r-$th receiver and $c(\bx)$ is the wave speed. The simulated domain $\Omega\subset\br^d$, $d$ is the dimension of the problem and $\bn$ is the unit outer normal vector to the boundary $\partial \Omega$. The point source hypothesis $\delta(\bx-\bxi)$ is considered in \eqref{eqn:wave} since we investigate the situation that the temporal and spatial scales of seismic wave propagated is large enough compared to the scales of seismic rupture \cite{AkRi:80,Ma:15}. The source time function has the form of Ricker wavelet
\begin{equation} \label{eqn:rick}
	f(t)=A\brac{1-2\pi^2f_0^2t^2}e^{-\pi^2f_0^2t^2},
\end{equation}
in which $f_0$ is the dominant frequency and $A$ is the normalization factor. Since we focus on the earthquake location problem on a large computational domain $\Omega$, we can simply consider the reflection boundary condition \eqref{bc:wave} rather than others, e.g. the perfectly matched layer absorbing boundary condition \cite{KoTr:03}.

The equations \eqref{eqn:optim}-\eqref{eqn:rick} provide a mathematical model to the waveform based earthquake location problem, which is a fundamental problem \cite{Th:14} with various applications in seismology \cite{LeSt:81,SaLoZo:08,WaEl:00}. Traditionally, the earthquake location problem is solved within the framework of the ray theory, see for example \cite{Ge:03,Ge:03b,Ge:12,PrGe:88,Th:85}. But the earthquake location results are not satisfactory since the ray theory is low accuracy when the seismic wave length is not small enough compared to the scale of wave propagation region \cite{EnRu:03,JiWuYa:08,RaPoFi:10,WuYa:13}. Thus, it is necessary to develop the waveform based earthquake location method. This direction is becoming more and more popular in recent years \cite{HuYaToBaQi:16,KiLiTr:11,LiGu:12,LiPoKoTr:04,ToZhYaYaChLi:14,WuChHuYa:16}, together with the fast increase of computational power. 

In general, the optimization problem \eqref{eqn:optim}-\eqref{eqn:wave} can be solved by iterative methods \cite{ToYaLiYaHa:16, WuChHuYa:16}. This approach requires that the initial value should be close enough to the global optimization solution. In particular, due to the highly singular nature of the delta function $\delta(\bx-\bxi)$ in the wave equation \eqref{eqn:wave}, the convergence domain could be very small. Although this issue has been studied in \cite{WuChHuYa:16}, the convergence domain is still restricted. The other important problem is the computational cost of the iterative method. To update the value of hypocenter and origin time, the sensitivity kernel should be obtained for each receiver. It costs several times wave equation computations to construct the sensitivity kernel. Considering the number of iterations and the number of initial values, the total computational cost could be very large. In practice, there are often two situations: the real-time earthquake location and many earthquakes relocation. Both of them require an accurate and efficient algorithm, which is the goal of this paper.

In this study, we introduce a set of auxiliary functions $\Xi_r(\bzeta,\nu)$. Define the solution set of the system
\begin{equation*}
	\mct=\{(\bzeta,\nu)\;|\; \Xi_r(\bzeta,\nu)=0,\;\forall r\in\mcr\},
\end{equation*}
in which $\mcr\subseteq \mca$, and $\mca$ is the set of all receivers. We will prove that $(\bxi_T,\tau_T)\in\mct$, i.e.
\begin{equation*}
	\Xi_r(\bxi_T,\tau_T)=0,\quad\forall r\in\mcr.
\end{equation*}
Therefore, the optimization problem \eqref{eqn:optim} is transformed into a system of equations, which can be solved in the least square sense. In the following sections, we will show that the computational cost of the AFM is significantly less than that of the iterative methods, e.g. \cite{LiPoKoTr:04,ToYaLiYaHa:16,WuChHuYa:16}.

The paper is organized as follows. In Section \ref{sec:taf}, we prove the main theorem and propose the algorithms. The numerical  experiments are presented to illustrate the features and highlights of AFM in Section \ref{sec:num}. In Section \ref{sec:con}, we make some conclusive remarks.

\section{The auxiliary function method} \label{sec:taf}

\subsection{The main theorem}
We begin this section with the following theorem. In fact, all the discussions in this paper depends on it.
\begin{theorem}
	For any given initial hypocenter $\bxi$ and origin time $\tau$, define the auxiliary functions
	\begin{equation} \label{fun:aux}
		\Xi_r(\bzeta,\nu)=2\chi_r(\bxi,\tau)-\int_0^Tf(t-\nu)w_r(\bzeta,t)-f(t-\tau)w_r(\bxi,t)\rd t, \quad \forall r\in\mcr,
	\end{equation}
	in which $\chi_r(\bxi,\tau)$ has been given in \eqref{fun:mis} and $w_r(\bx,t)$ satisfies the adjoint equation with terminal-boundary conditions
	\begin{equation} \label{eqn:adjoint}
		\left\{\begin{array}{ll}
			\frac{\partial^2 w_r(\bx,t)}{\partial t^2}=\nabla\cdot\brac{c^2(\bx)\nabla w_r(\bx,t)}
				+\frac{d_r(t)-s(\boeta_r,t)}{\int_0^T\abs{d_r(t)}^2\rd t}\delta(\bx-\boeta_r), & \bx\in\Omega, \\
	                 w_r(\bx,T)=\frac{\partial w_r(\bx,T)}{\partial t}=0, & \bx\in\Omega, \\
        		         	 \bn\cdot\brac{c^2(\bx)\nabla w_r(\bx,t)}=0, & \bx\in\partial\Omega.
		\end{array}\right.
	\end{equation}
	Then, for real earthquake hypocenter $\bxi_T$ and origin time $\tau_T$, we have
	\begin{equation}
		\Xi_r(\bxi_T,\tau_T)=0,\quad\forall r\in\mcr.
	\end{equation}	
\end{theorem}
{\bf Proof:} Let us first define the difference function
\begin{equation*}
	\delta s(\bx,t)=u(\bx,t;\bxi_T,\tau_T)-u(\bx,t;\bxi,\tau).
\end{equation*}
According to \eqref{eqn:rel_ds}, it follows
\begin{equation*}
	\delta s(\boeta_r,t)=d_r(t)-s(\boeta_r,t), \quad \forall r\in\mcr.
\end{equation*}
And the difference function $\delta s(\bx,t)$ satisfies the wave equation with initial-boundary conditions
\begin{equation}  \label{eqn:diff_fun}
	\left\{\begin{array}{ll}
		\frac{\partial^2 \delta s(\bx,t)}{\partial t^2}=\nabla\cdot\brac{c^2(\bx)\nabla \delta s(\bx,t)}
			+f(t-\tau_T)\delta(\bx-\bxi_T)-f(t-\tau)\delta(\bx-\bxi), & \bx\in\Omega, \\
                \delta s(\bx,0)=\frac{\partial \delta s(\bx,0)}{\partial t}=0, & \bx\in\Omega, \\
                \bn\cdot\brac{c^2(\bx)\nabla \delta s(\bx,t)}=0, & \bx\in\partial\Omega.
        \end{array}\right.
\end{equation}
Multiply the wave function $w_r(\bx,t)$ given in \eqref{eqn:adjoint}, integrate it on $\Omega\times[0,T]$ and use the integration by parts, we obtain
\begin{multline} \label{eqn:integration}
	\int_0^T\int_{\Omega}\frac{\partial^2 w_r}{\partial t^2}\delta s\rd \bx\rd t
		=\int_0^T\int_{\Omega}\delta s\nabla\cdot(c^2\nabla w_r)\rd \bx \rd t \\
                 +\int_0^T f(t-\tau_T)w_r(\bxi_T,t)-f(t-\tau)w_r(\bxi,t)\rd t.
\end{multline}
On the other hand, the misfit function $\chi_r(\bxi,\tau)$ in \eqref{fun:mis} can be rewritten as
\begin{multline} \label{fun:mis2}
	\chi_r(\bxi,\tau)=\frac{\int_0^T\brac{d_r(t)-s(\boeta_r,t)}\delta s(\boeta_r,t)\rd t}{2\int_0^T\abs{d_r(t)}^2\rd t} \\
		=\frac{\int_0^T\int_{\Omega}\brac{d_r(t)-s(\boeta_r,t)}\delta s(\bx,t)\delta(\bx-\boeta_r)\rd\bx\rd t}{2\int_0^T\abs{d_r(t)}^2\rd t}.
\end{multline}
Multiplying both sides of the above equation by $2$, and adding equation \eqref{eqn:integration}, we get
\begin{equation*}
	2\chi_r(\bxi,\tau)=\int_0^T f(t-\tau_T)w_r(\bxi_T,t)-f(t-\tau)w_r(\bxi,t)\rd t.
\end{equation*}
This completes the proof. $\Box$

This theorem holds for any initial hypocenter $\bxi$ and origin time $\tau$. Assume that the solution of the equations
\begin{equation} \label{eqn:Xi_root}
	\Xi_r(\bzeta,\nu)=0, \quad \forall r\in\mcr,
\end{equation}
is unique, we only need one round computation to get the real hypocenter $\bxi_T$ and origin time $\tau_T$. In this sense, the theorem leads to a non-iterative method. The theorem does not require that the real hypocenter $\bxi_T$ and the initial hypocenter $\bxi$ are close. Thus, a global method for the inverse problem \eqref{eqn:optim} can be developed based on the theorem.

There are two issues we need to clarify about the solutions of equation \eqref{eqn:Xi_root}:
\begin{enumerate}[1.]
	\item In practice, we prefer to solve the equations  \eqref{eqn:Xi_root} in the least square sense, i.e.
	\begin{equation} \label{eqn:Xi_ls}
		(\bxi_T,\tau_T)=\argmin_{\bzeta,\nu}\Gamma(\bzeta,\nu).
	\end{equation}
	in which
	\begin{equation*}
		\Gamma(\bzeta,\nu)=\sum_{r\in\mcr}\Xi_r^2(\bzeta,\nu).
	\end{equation*}

	\item The theorem does not guarantee the uniqueness of the solution. This may lead to incorrect inversion result. Fortunately, more constraints $r\in\mcr$ may improve the uniqueness of the solution. According to the numerical experiments, we are not suffering from the problem of uniqueness. $\Box$
\end{enumerate}

\subsection{Algorithm and Discussions}
Based on the above results, the detailed algorithm is as follows:

\begin{algorithm}[The auxiliary function method] \label{alg:afm_class} $ $
\begin{itemize}
	\item Initialization. Given a searching domain $\Omega_s$ and a searching time interval $I_s$, we wish that $\bxi_T\in\Omega_s$ and $\tau_T\in I_s$. Select a mesh size $h>0$ and a time step $\sigma>0$.
	
	\item Discretization. Select a uniform or quasi-uniform grid $\Omega_h\subset \Omega_s$ of mesh size $h$. Select a uniform or quasi-uniform time division $I_{\sigma}\subset I_s$ of time step $\sigma$.
	
	\item Forward and Backward evolution. Given the initial hypocenter $\bxi$ and origin time $\tau$, compute the wave equation \eqref{eqn:wave}-\eqref{bc:wave} to get
		\begin{equation*}
			s(\boeta_r,t)=u(\boeta_r,t;\bxi,\tau),
		\end{equation*}
		and the misfit function $\chi_r(\bxi,\tau),\;\forall r\in\mcr$. Next, compute the individual adjoint wave equation \eqref{eqn:adjoint} to get $w_r(\bx,t)$ for $r\in\mcr$. 
	
	\item Construction. Evaluate all the values of the auxiliary functions $\Xi_r(\bzeta,\nu)$ on the mesh size $\bzeta\in\Omega_h$, time division $\nu\in I_{\sigma}$ and $r\in\mcr$. Thus, $\Gamma(\bzeta,\nu)$ for $(\bzeta,\nu) \in \Omega_h \times I_{\sigma}$ is directly obtained.
	
	\item Output. Finally, we can easily get the approximated solution of the optimization problem \eqref{eqn:Xi_ls}
	\begin{equation*}
		(\bxi_*,\tau_*)=\argmin_{\bzeta\in\Omega_h,\nu\in I_{\sigma}}\Gamma(\bzeta,\nu).
	\end{equation*}
	by direct search. Output $(\bxi_*,\tau_*)$ and stop. $\Box$
\end{itemize}
\end{algorithm}

Now, we get an estimate of the real hypocenter and origin time $(\bxi_*,\tau_*)$. But we still need to check whether this is valid. For the case where 
\begin{equation*}
	\bxi_T\notin\Omega_s\quad\textrm{or}\quad \tau_T\notin I_s,
\end{equation*}
it is obviously that the output $(\bxi_*,\tau_*)$ is a wrong approximation. To avoid this situation, we have to check that if
\begin{equation} \label{eqn:criteria}
	\sum_{r\in\mcr}\abs{\chi_r(\bxi_*,\tau_*)}<\veps_1,
\end{equation}
is satisfied. Here $\veps_1$ is the tolerance value. To obtain the values of the misfit function $\chi_r(\bxi_*,\tau_*)$ for all $r\in\mcr$, we only need to solve the forward wave equation for one time. Once the criteria \eqref{eqn:criteria} is met, an accurate approximation of the real hypocenter and origin time is obtained. Otherwise, we need to restart the algorithm with different searching domain $\Omega_s$ and searching time interval $I_s$.

\begin{remark} \label{rem:conv_cond}
	In fact, it is not expensive to satisfy the conditions
	\begin{equation*}
		\bxi_T\in\Omega_s\quad\textrm{and}\quad \tau_T\in I_s.
	\end{equation*}
	This only needs a large searching domain $\Omega_s$ and searching time interval $I_s$. According to the above algorithm, we don't need to solve extra wave equations for this. $\Box$
\end{remark}

We now discuss the computational cost of Algorithm \ref{alg:afm_class}, which consists of four parts:
\begin{enumerate}[1.]
	\item The simulation of the wave equation \eqref{eqn:wave}-\eqref{bc:wave} and the adjoint wave equations \eqref{eqn:adjoint}. It needs to solve the wave equation for $\#\mcr+1$ times. Here $\#\mcr$ denotes the number of elements in the set $\mcr$.
	
	\item The computation of the misfit functions $\chi_r(\bxi,\tau)$. It needs to compute the $1$-d integral for $2\#\mcr$ times.
	
	\item The computation of the auxiliary functions $\Xi_r(\bzeta,\nu)$. It needs to compute the $1$-d integral for $\#\mcr\cdot(\#\Omega_h\cdot\#I_{\sigma}+1)$ times.
	
	\item Computing the misfit function $\chi_r(\bxi_*,\tau_*)$ for all $r\in\mcr$. It requires to solve the forward wave equation once and calculate the $1$-d integral $2\#\mcr$ times. $\Box$
\end{enumerate}

Adding these four parts, the total computational cost is $\#\mcr+2$ times wave equation computations plus $\#\mcr\cdot(\#\Omega_h\cdot\#I_{\sigma}+5)$ times $1$-d integral calculations. Considering the following facts: (i) the support of integrand function is small, (ii) the mesh size and the time step of the searching domain and the searching time interval can be coarser than that of the wave field simulations, the overall cost does not exceed $\#\mcr+3$ times wave equation computations. 

For iterative method, the overall computational cost is about $(\#\mcr+1)\cdot N_{iter}\cdot N_{init}$ times wave equation computations. Here $N_{iter}$ denotes the average number of iterations for each initial data and $N_{init}$ denotes the average number of initial data for each problem. Obviously, the AFM is much more efficient.

\subsection{The auxiliary function preprocessor method}

There is one problem with Algorithm \ref{alg:afm_class}. In order to improve the location accuracy, we need to reduce the mesh size $h$ and time step $\sigma$. This will increase the computational cost of the Construction step. In the extreme case, the computational cost of the Construction step could be far more than the cost of wave equation computations. It is of course not worthwhile. An alternative way to improve the location accuracy is to treat the auxiliary function method (Algorithm \ref{alg:afm_class}) as a preprocessor. The sketch of the algorithm is as follows:

\begin{algorithm}[The auxiliary function preprocessor method] \label{alg:afm_preproc} $ $
\begin{itemize}
	\item Initialization. Given a searching domain $\Omega_s$ and a searching time interval $I_s$.

	\item Preprocessing. Execute Algorithm \ref{alg:afm_class} on a coarse mesh grid $\Omega_h\subset \Omega_s$ and time division $I_{\sigma}\subset I_s$ to obtain $(\bxi_*,\tau_*)$.
	
	\item Verifying. If the criterion \eqref{eqn:criteria} is not met, go back to the Preprocessing step with different searching domain $\Omega_s$ and searching time interval $I_s$. Otherwise, go to the Iteration step.
	
	\item Iteration. Execute the iterative method, see e.g. in \cite{ToYaLiYaHa:16,WuChHuYa:16}, for the inverse problem \eqref{eqn:optim}-\eqref{eqn:rick} with initial hypocenter and origin time $(\bxi_*,\tau_*)$ to get $(\bxi_A,\tau_A)$. $\Box$
\end{itemize}
\end{algorithm}

According to the above algorithm, we get an accurate estimation of the earthquake hypocenter and origin time $(\bxi_A,\tau_A)$. The total computational cost of this algorithm is about $(\#\mcr+3)\cdot N_{AFM}+(\#\mcr+1)\cdot N_{iter}'$ times wave equation computations. Here $N_{AFM}$ and $N_{iter}'$ denote the number of execution times for AFM and the average number of iterations of the iterative method. As discussed in Remark \ref{rem:conv_cond}, it is not difficult to satisfy the criterion \eqref{eqn:criteria}. Thus $N_{AFM}$ can be considered as $1$. Moreover, since the initial value $(\bxi_*,\tau_*)$ is close enough to optimization solution $(\bxi_T,\tau_T)$, the number of iterations $N_{iter}'$ can also be very small. Thus, the AFM preprocessor is still more efficient than the iterative methods.

\section{Numerical Experiments} \label{sec:num}
In this section, two examples are presented to demonstrate the efficiency and effectiveness of our method. In all the numerical examples, the finite difference schemes \cite{Da:86,YaLuWuPe:04} are applied to solve the acoustic wave equation \eqref{eqn:wave} with initial condition \eqref{ic:wave}. On the surface of the earth, we consider the reflection boundary condition \eqref{bc:wave}. And the perfectly matched layer boundary condition \cite{KoTr:03} is used for the other boundaries within the earth. The point source $\delta(x-\xi)$ is numerical discretized as follows \cite{We:08}:
\begin{equation*}
	\delta(x)=\left\{\begin{array}{ll}
		\frac{1}{h}\left(1-\frac{5}{4}\abs{\frac{x}{h}}^2-\frac{35}{12}\abs{\frac{x}{h}}^3
			+\frac{21}{4}\abs{\frac{x}{h}}^4-\frac{25}{12}\abs{\frac{x}{h}}^5\right), & \abs{x}\le h, \\
		\frac{1}{h}\left(-4+\frac{75}{4}\abs{\frac{x}{h}}-\frac{245}{8}\abs{\frac{x}{h}}^2+\frac{545}{24}\abs{\frac{x}{h}}^3
			-\frac{63}{8}\abs{\frac{x}{h}}^4+\frac{25}{24}\abs{\frac{x}{h}}^5\right), & h<\abs{x}\le 2h, \\
		\frac{1}{h}\left(18-\frac{153}{4}\abs{\frac{x}{h}}+\frac{255}{8}\abs{\frac{x}{h}}^2-\frac{313}{24}\abs{\frac{x}{h}}^3
			+\frac{21}{8}\abs{\frac{x}{h}}^4-\frac{5}{24}\abs{\frac{x}{h}}^5\right), & 2h<\abs{x}\le 3h, \\
		0, & \abs{x}>3h.
	\end{array}\right. 
\end{equation*}
Here $h$ is a numerical parameter which is related to the mesh size.

\subsection{The two-layer velocity model} \label{subsec:two_layer}
\begin{figure} 
	\centering
	\includegraphics[width=0.60\textwidth, height=0.16\textheight]{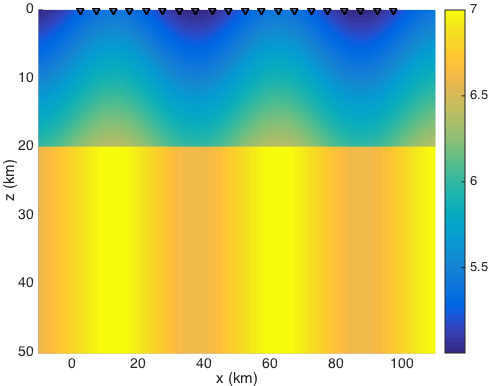}
	\caption{Illustration of two-layer velocity model. The black triangles indicate the receivers.} \label{fig:exam31_vel}
\end{figure}

Consider the two-layer model in the bounded domain $\Omega=[-10\,km,110\,km]\times[0\,km, 50\,km]$, the wave speed is
\begin{equation*}
	c(x,z)=\left\{\begin{array}{ll}
		5.2+0.05z+0.2\sin\frac{\pi x}{25}, & 0\,km\le z\le 20\,km, \\
		6.8+0.2\sin\frac{\pi x}{25}, & z>20\,km.
	\end{array}\right.
\end{equation*}
The unit is `km/s'. The computational time interval $I=[0,25\,s]$. The dominant frequency of the earthquakes is $f_0=2\,Hz$. There are $20$ equidistant receivers on the surface
\begin{equation*}
	\boeta_r=(x_r,z_r)=(5r-2.5\,km,0), \quad r=1,2,\cdots,20,
\end{equation*}
see Figure \ref{fig:exam31_vel} for illustration.

First, we test the auxiliary function preprocessor method (Algorithm \ref{alg:afm_preproc}) using $500$ experiments. The searching domain is $\Omega_s=[0,100\,km]\times[0,40\,km]$, and the searching time interval is $I_s=[0,25\,s]$. The mesh sizes for the searching grid are $h_x=0.5\,km$ and $h_z=0.4\,km$. The time step for the searching time interval is $\sigma=0.1\,s$. As a comparison, we also compute these experiments by the iterative method proposed in \cite{WuChHuYa:16}.

\begin{figure} 
	\centering
	\includegraphics[width=0.90\textwidth, height=0.16\textheight]{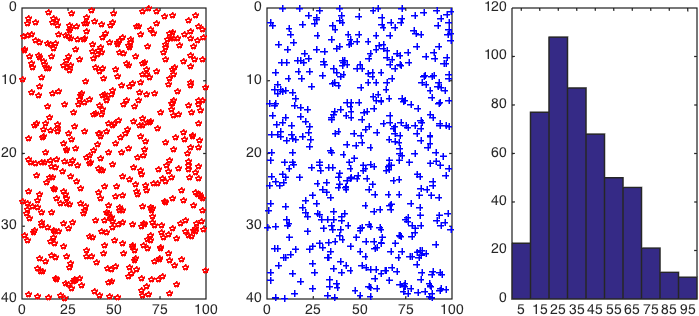}
	\caption{The two-layer velocity model. Left: the spatial distribution of the real earthquake hypocenter $\bxi_T^i$; Middle: the spatial distribution of the initial earthquake hypocenter $\bxi^i$; Right: the distance distribution histogram between the real and the initial earthquake hypocenter $d^i$.} \label{fig:exam31_dis}
\end{figure}

The experiments are designed as follows: the real and initial earthquake hypocenter $\bxi_T^i,\;\bxi^i$ are both uniformly distributed over $[0,100\,km]\times[0,40\,km]$, the real and initial original time $\tau_T^i,\;\tau^i$ are both uniformly distributed over $[5\,s,20\,s]$. Their spatial distribution and the histogram of the distance between the real and the initial hypocenter
\begin{equation*}
	d^i=\norm{\bxi_T^i-\bxi^i}_2,
\end{equation*}
 are presented in Fig \ref{fig:exam31_dis}.

In the auxiliary function preprocessor method, we randomly select five receivers for inversion, e.g. $r=3,\;5,\;9,\;14,\;18$. In Table \ref{tab:exam31_convergent}, we can see the convergent results of the two methods. From which, we can conclude that the auxiliary function preprocessor method converges globally here. For the iterative method proposed in \cite{WuChHuYa:16}, only $23.4\%$ experiments converge. We have to remark that the range of convergence of the iterative method has been enlarged by several tens of time. In contrast, the auxiliary function preprocessor method has an absolute advantage in terms of convergence.

\begin{table*}
    	\caption{The two-layer velocity model. Convergent results for the auxiliary function preprocessor method(AFPM) and the iterative method(IM).} \label{tab:exam31_convergent}
	\begin{center}\begin{tabular}{ccccc} \hline
		 & Correct convergence & Diverge & Error convergence & Total \\ \hline
		 AFPM & $500$ & $0$ & $0$ & $500$\\ \hline
		 IM & $117$ & $355$ & $28$ & $500$ \\ \hline
    	\end{tabular}\end{center}
\end{table*}

\begin{figure} 
	\centering
	\includegraphics[width=0.90\textwidth, height=0.30\textheight]{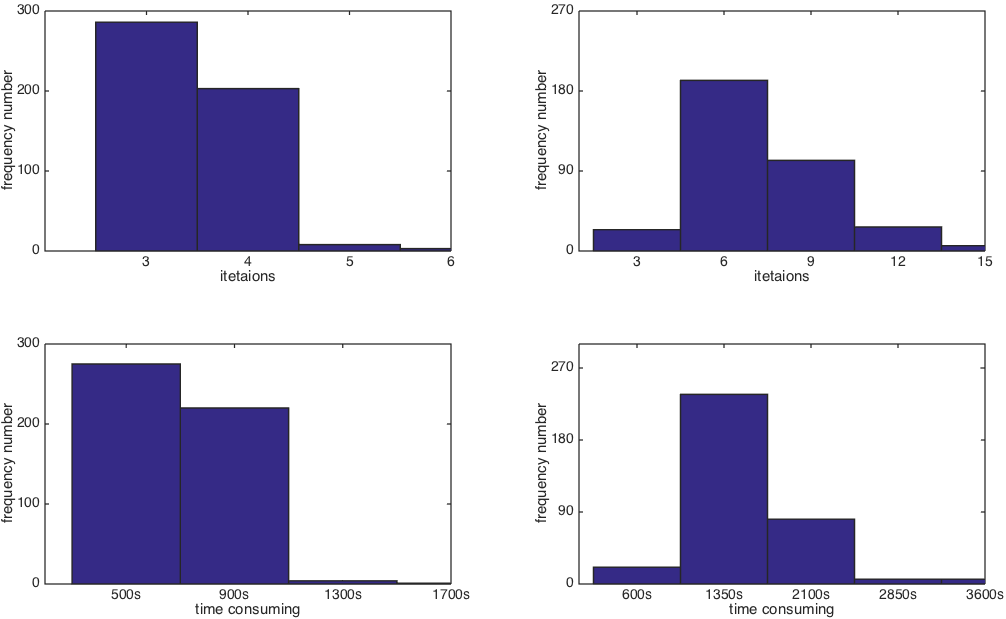}
	\caption{The two-layer velocity model. The histogram of the iterations and computational time for the two methods. Up: the histogram of iterations; Down: the histogram of computational time; Left: the auxiliary preprocessor method; Right: the iterative method.} \label{fig:exam31_itti}
\end{figure}

In Figure \ref{fig:exam31_itti}, we output the histogram of the iterations and computational time for the two methods. The mean and standard deviation of iterations and computational time for the two methods are also presented in Table \ref{tab:exam31_iterationandtime}. In particular, we present the mean and standard deviation of the time consuming for the preprocessing step and the single iteration step of the auxiliary function preprocessor method in Table \ref{tab:exam31_meanvariance}. It is obviously that the time consuming of the preprocessing step and the single iteration step are almost the same. Thus, we consider the preprocessing step to be an iteration step. Taking account of all the above issues, the total computational cost of the iterative method is about
\begin{equation*}
	\frac{500}{117}\times \frac{1523}{689}\approx 9.45 \; \textrm{times}
\end{equation*}
of the auxiliary function preprocessor method. Thus, we can concluded that the auxiliary function preprocessor method is more efficient than the iterative method. This agrees with the theoretical discussions in the previous section.

\begin{table*}
    	\caption{The two-layer velocity model. The Mean(M) and Standard Deviation(SD) of iterations and computational time for the AFPM and IM.} \label{tab:exam31_iterationandtime}
	\begin{center}\begin{tabular}{c|cc|cc} \hline
		 & \multicolumn{2}{c|}{Iterations} & \multicolumn{2}{c}{Computational time} \\
		 & M & SD & M & SD \\ \hline
		 AFPM & $3.46$ & $0.61$ & $689\,s$ & $132.7\,s$ \\ \hline
		 IM & $8.33$ & $2.35$ & $1523\,s$ & $501.0\,s$ \\ \hline
    	\end{tabular}\end{center}
\end{table*}

\begin{table*}
    	\caption{The two-layer velocity model. The mean and standard deviation of the time consuming for the preprocessing step and the single iteration step of the auxiliary function preprocessor method.} \label{tab:exam31_meanvariance}
	\begin{center}\begin{tabular}{ccc} \hline
		 & Mean & Standard Deviation \\ \hline
		 The preprocessing step & $224.9\,s$ & $13.2\,s$ \\ \hline
		 The single iteration step & $222.7\,s$ & $11.8\,s$ \\ \hline
    	\end{tabular}\end{center}
\end{table*}

\begin{figure*} 
	\begin{tabular}{ccc}
		\includegraphics[width=0.33\textwidth, height=0.12\textheight]{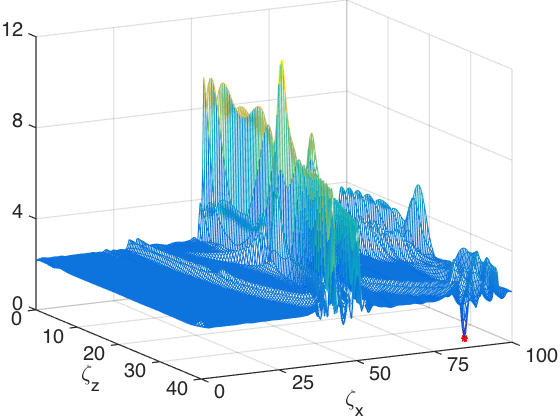} &
		\includegraphics[width=0.33\textwidth, height=0.12\textheight]{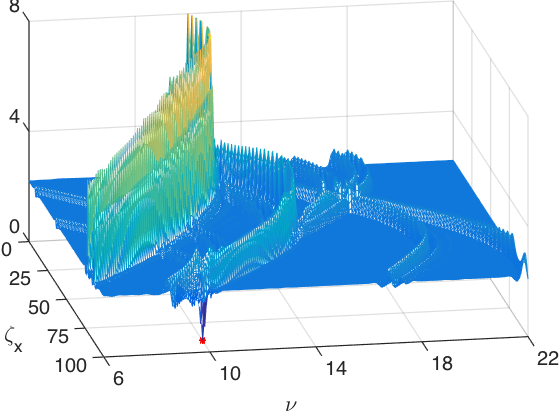} &
		\includegraphics[width=0.33\textwidth, height=0.12\textheight]{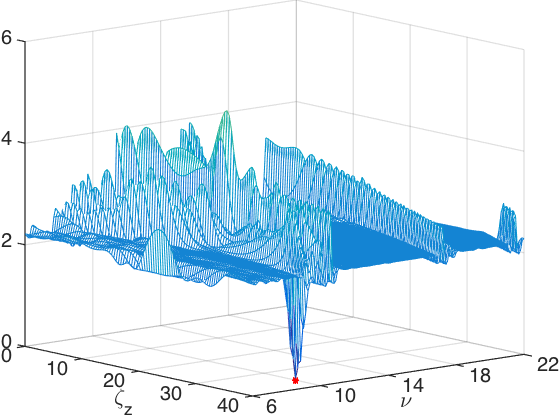} \\
		\includegraphics[width=0.33\textwidth, height=0.12\textheight]{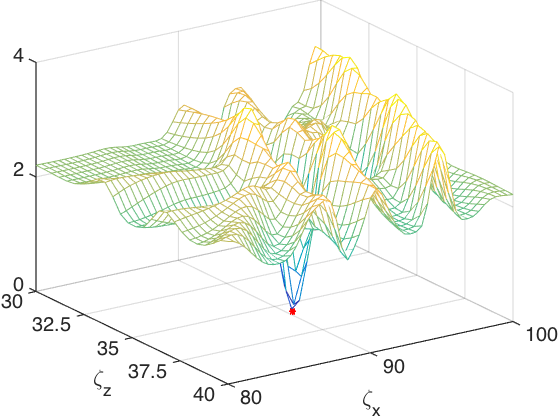} &
		\includegraphics[width=0.33\textwidth, height=0.12\textheight]{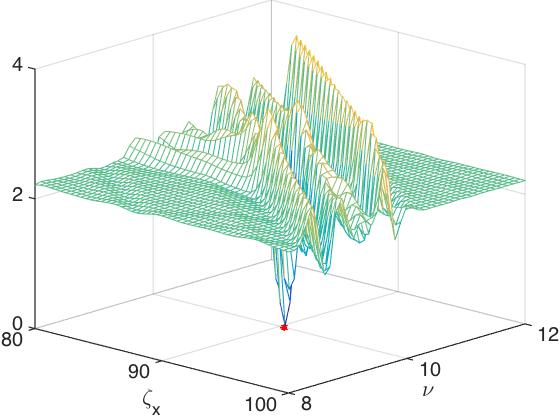} &
		\includegraphics[width=0.33\textwidth, height=0.12\textheight]{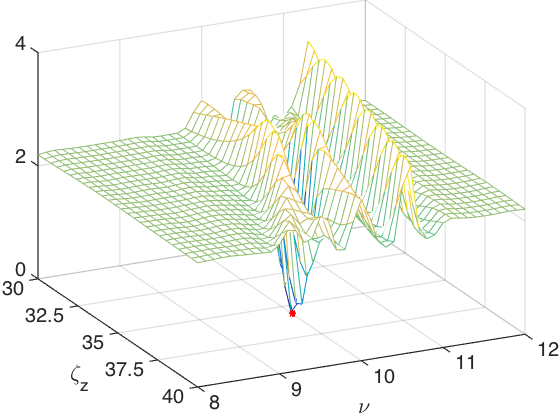}
	\end{tabular}
	\caption{The two-layer velocity model, case (i). The $\zeta_x-\zeta_z$ (Left, $\nu=10\,s$), $\nu-\zeta_x$ (Middle, $\;\zeta_z=35.67\,km$) and $\nu-\zeta_z$ (Right, $\zeta_x=90.36\,km$) cross section plan of function $\Gamma(\bzeta,\nu)$. Up for large image and Down for zoom-in image.} \label{fig:exam31_2exam_c1auxfun}
\end{figure*}

\begin{figure*} 
	\begin{tabular}{ccc}
		\includegraphics[width=0.33\textwidth, height=0.12\textheight]{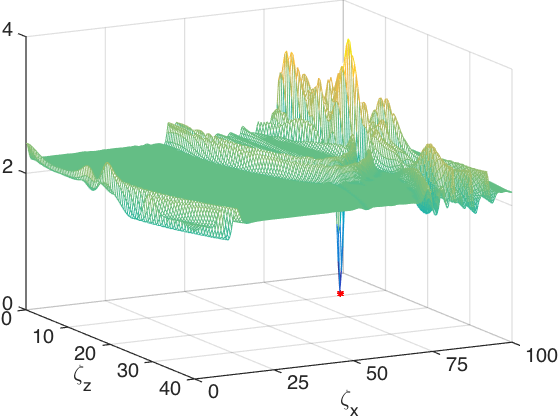} &
		\includegraphics[width=0.33\textwidth, height=0.12\textheight]{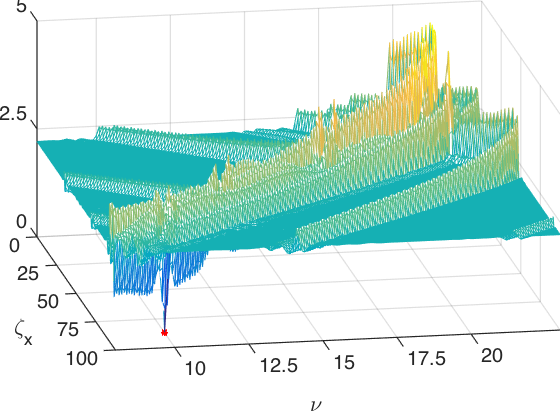} &
		\includegraphics[width=0.33\textwidth, height=0.12\textheight]{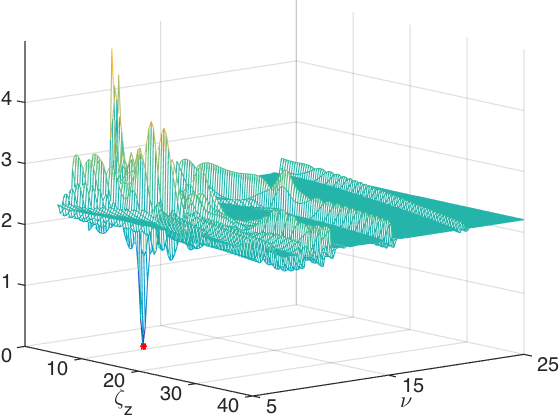} \\
		\includegraphics[width=0.33\textwidth, height=0.12\textheight]{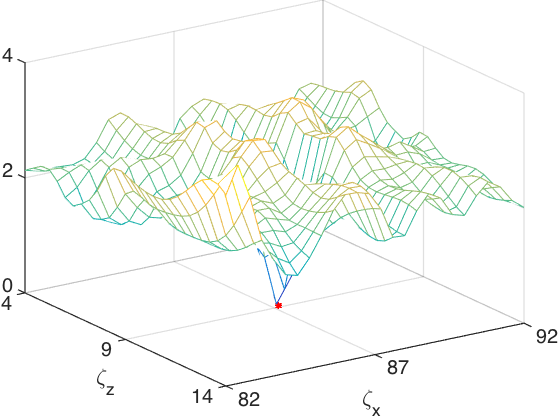} &
		\includegraphics[width=0.33\textwidth, height=0.12\textheight]{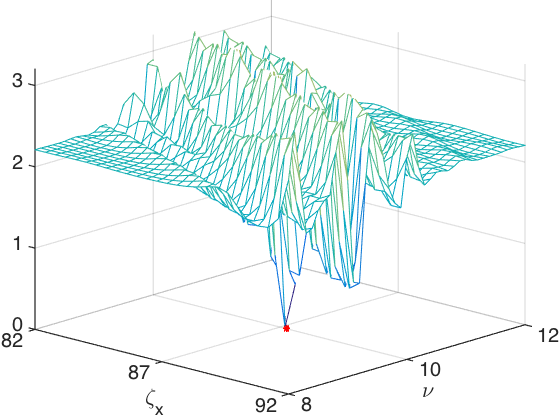} &
		\includegraphics[width=0.33\textwidth, height=0.12\textheight]{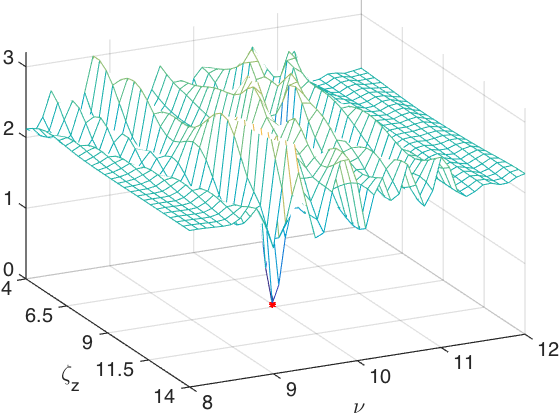}
	\end{tabular}
	\caption{The two-layer velocity model, case (ii). The $\zeta_x-\zeta_z$ (Left, $\nu=10\,s$), $\nu-\zeta_x$ (Middle, $\;\zeta_z=8.842\,km$) and $\nu-\zeta_z$ (Right, $\zeta_x=87.252\,km$) cross section plan of function $\Gamma(\bzeta,\nu)$. Up for large image and Down for zoom-in image.} \label{fig:exam31_2exam_c2auxfun}
\end{figure*}

\begin{figure*} 
	\begin{tabular}{c}
		\includegraphics[width=0.95\textwidth, height=0.12\textheight]{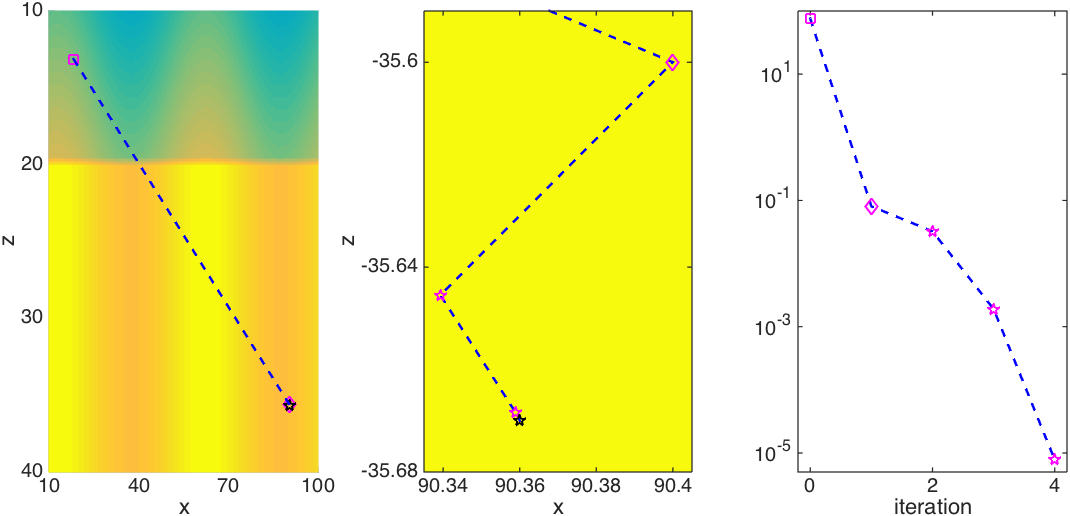} \\
		\includegraphics[width=0.95\textwidth, height=0.12\textheight]{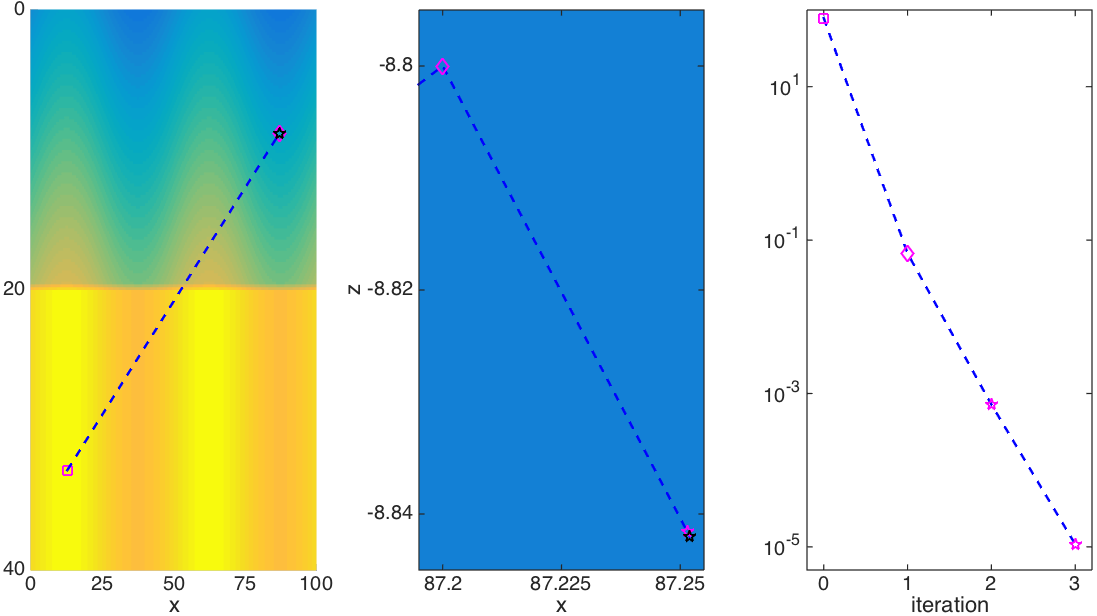}
	\end{tabular}
	\caption{Convergence history of the two-layer velocity model.  Up for case (i), and Down for case (ii). Left (large image) and Middle (zoom-in image): the convergent trajectories; Right: the absolute errors with respect to iteration steps between the real and computed hypocenter of the earthquake. The magenta square is the initial hypocenter, the magenta diamond denotes the hypocenter obtained by AFM (algorithm \ref{alg:afm_class}), the magenta pentagrams indicate the hypocenters in the iterative process of the AFPM (algorithm \ref{alg:afm_preproc}), and the black pentagram is the real hypocenter.} \label{fig:exam31_2exam_convtraj}
\end{figure*}

Next, two examples are specifically presented. The parameters are selected as follows:
\begin{align*}
	(i) \; \bxi_T=(90.36\,km,35.67\,km),\; \tau_T=10\,s, \quad \bxi=(18.23\,km,13.13\,km), \; \tau=15.5\,s; \\
	(ii) \; \bxi_T=(87.252\,km,8.842\,km),\; \tau_T=10\,s, \quad \bxi=(12.75\,km,32.87\,km), \; \tau=17.4\,s.
\end{align*}
In Figure \ref{fig:exam31_2exam_c1auxfun}-\ref{fig:exam31_2exam_c2auxfun}, we output the functions $\Gamma(\bzeta,\nu)$, from which we can see the global minimum is unique in both cases. The convergent history of the auxiliary preprocessor method are also illustrated in Figure \ref{fig:exam31_2exam_convtraj}. Here, we randomly select five receivers for inversion, e.g. $r=3,\;5,\;9,\;14,\;18$. We can see that the global minimum of the function $\Gamma(\bzeta,\nu)$ is very close to the optimization solution of \eqref{eqn:optim}. Thus, when the accuracy requirement is not high, the solution of auxiliary function method is directly applicable. On the other hand, when the accuracy requirement is high, the solution of auxiliary function method can provide excellent initial values for the iterative methods.

\begin{figure*} 
	\begin{tabular}{cccc}
		\includegraphics[width=0.23\textwidth, height=0.12\textheight]{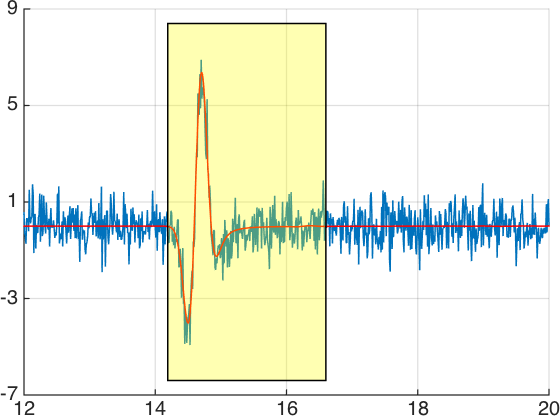} &
		\includegraphics[width=0.23\textwidth, height=0.12\textheight]{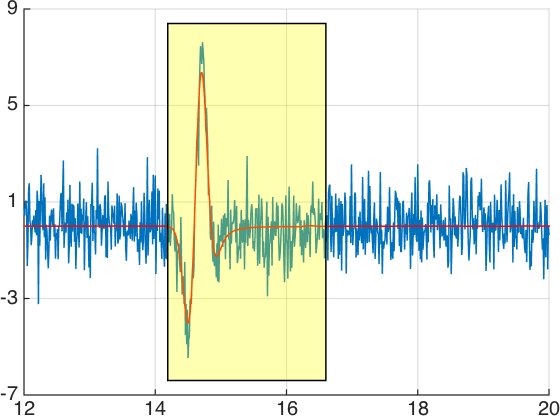} &
		\includegraphics[width=0.23\textwidth, height=0.12\textheight]{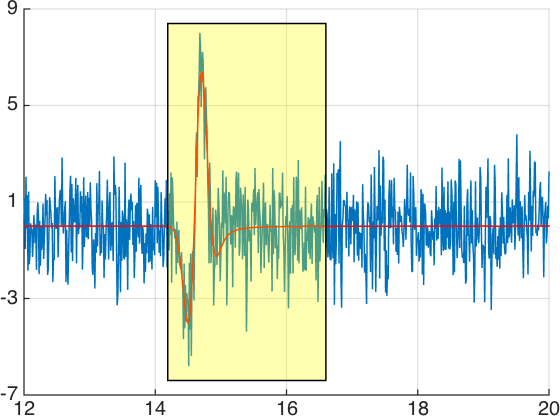} &
		\includegraphics[width=0.23\textwidth, height=0.12\textheight]{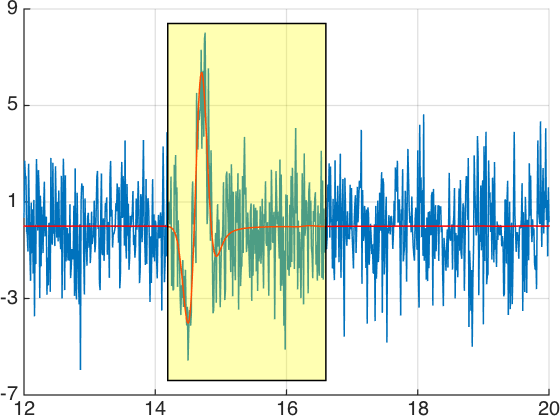} \\
		\includegraphics[width=0.23\textwidth, height=0.12\textheight]{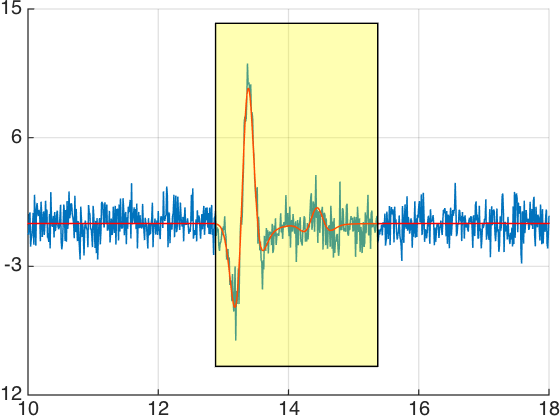} &
		\includegraphics[width=0.23\textwidth, height=0.12\textheight]{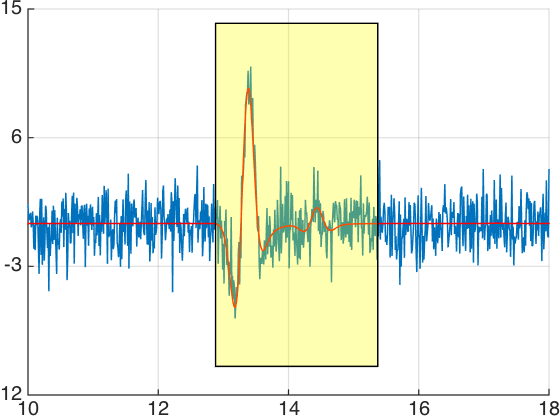} &
		\includegraphics[width=0.23\textwidth, height=0.12\textheight]{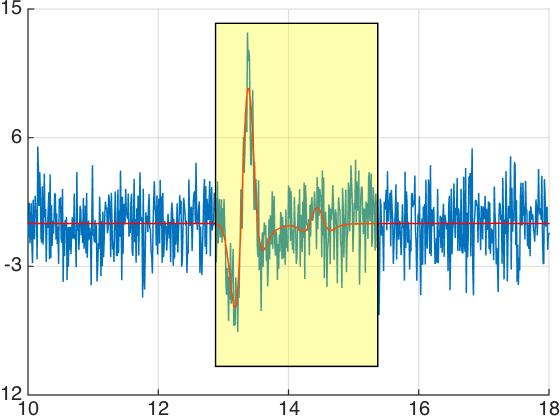} &
		\includegraphics[width=0.23\textwidth, height=0.12\textheight]{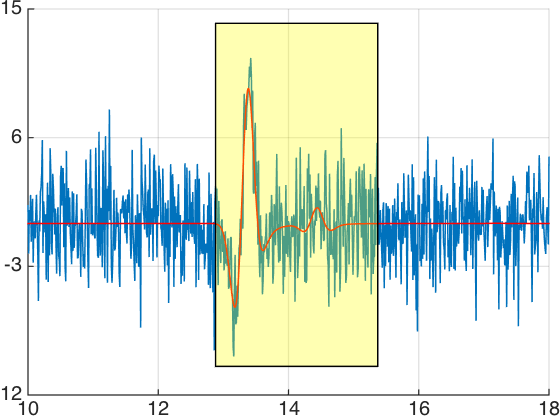}
	\end{tabular}
	\caption{Illustration of signal with noise in the two-layer velocity model. The signal with noise $d_r(t)$ (blue line) and the noise free signal $u(\boeta_r,t;\bxi_T,\tau_T)$ for receivers $r=9$. The horizontal axis is the time $t$. Up: parameters group (i); Down: parameters group (ii); From left to right, the ratio $R=10\%,\;15\%,\;20\%,\;25\%$ respectively.} \label{fig:exam31_noise_signal}
\end{figure*}

At last, we test the influence of the noise. The same parameters (i) and (ii) are selected here. The real earthquake signal can be regarded as
\begin{equation*}
	d_r(t)=u(\boeta_r,t;\bxi_T,\tau_T)+N_r(t),
\end{equation*}
with $N_r(t)$ is subject to the normal distribution with mean $\mu=0$ and the standard deviation
\begin{equation*}
	\sigma=R \times \max_t\abs{u(\boeta_r,t;\bxi_T,\tau_T)}.
\end{equation*}
Here $R$ denotes the ratio, which will be selected as $10\%,\;15\%,\;20\%$ and $25\%$ respectively. The real earthquake signal with noise $d_r(t)$ and the noise free signal $u(\boeta_r,t;\bxi_T,\tau_T)$ are illustrated in Fig \ref{fig:exam31_noise_signal}. In order to reduce the impact of noise, we can select a time window that contains the main part of $u(\boeta_r,t;\bxi_T,\tau_T)$. In Table \ref{tab:exam31_datanoise_error}, the mean and standard deviation of the errors between the location results $(\bxi_*,\tau_*)$ computed via the auxiliary function method (Algorithm \ref{alg:afm_class}) and the true solution $(\bxi_T,\tau_T)$ are presented. The auxiliary function preprocessor method (Algorithm \ref{alg:afm_preproc}) is not considered here since the iterative method may fail even for small ratio $R$. For each parameter group and ratio $R$, we test the algorithm with 10 different noise. We note that all the standard deviations are zero. This implies that all the test converge to the same solution. They are
\begin{align*}
	(i) \; \bxi_*=(90.40\,km,35.60\,km),\; \tau_*=10.01\,s; \\
	(ii) \; \bxi_*=(87.20\,km,8.80\,km),\; \tau_*=9.99\,s.
\end{align*}
This is because we are using the same mesh grid $\Omega_h$ and time division $I_{\sigma}$ in the program. It should be noted that the numerical results may become better or worse as the mesh grid and time division change. But the location errors are always in the same order of magnitude, which is caused by the noise. We can also observe that the algorithm failed when $R=25\%$ in the parameter group (i). But the algorithm works for all the ratios in the parameter group (ii). Thus, we tend to believe that the algorithm can be success when $R$ below $20\%$ in this example. This is a huge advantage of the auxiliary function method over the traditional iterative methods. In fact, accordingly to our numerical tests, the traditional iterative methods fail when $R$ reaches to $10\%$ and above. Moreover, the computational cost of the auxiliary function method is much less than that of the conventional iterative methods.

\begin{table*}
    	\caption{The two-layer velocity model. The Mean(M) and Standard Deviation(SD) of of the errors between the location results $(\bxi_*,\tau_*)$ and the true solution $(\bxi_T,\tau_T)$.} \label{tab:exam31_datanoise_error}
	\begin{center}\begin{tabular}{c|cc|cc} \hline
		 & \multicolumn{2}{c|}{Errors of AFM, case (i)} & \multicolumn{2}{c}{Errors of AFM, case (ii)} \\
		 $R$ & M & SD & M & SD \\ \hline
		 $10\%$ & $0.0812$ & $0$ & $0.0676$ & $0$ \\
		 $15\%$ & $0.0812$ & $0$ & $0.0676$ & $0$ \\
		 $20\%$ & $0.0812$ & $0$ & $0.0676$ & $0$ \\
		 $25\%$ & \multicolumn{2}{c|}{fail} & $0.0676$ & $0$ \\ \hline
    	\end{tabular}\end{center}
\end{table*}

\subsection{The practical velocity model} \label{subsec:practical}
\begin{figure} 
	\centering
	\includegraphics[width=0.60\textwidth, height=0.25\textheight]{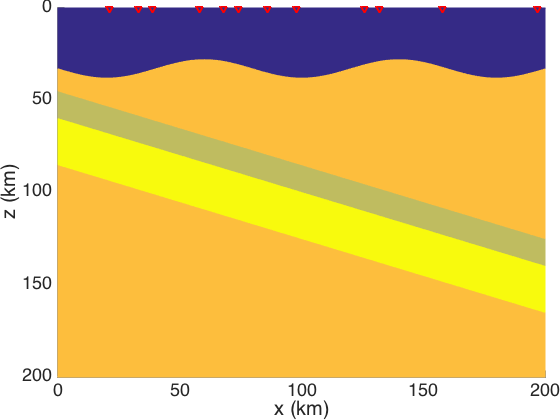}
	\caption{The practical velocity model. The red triangles indicate the receivers.} \label{fig:pracvel_illu}
\end{figure}

Let's consider a more practical model, the computational domain is $[0\,km,200\,km]\times[0\,km,200\,km]$, and the wave speed is 
	\begin{equation*}
		c(x,z)=\left\{\begin{array}{ll}
			5.5, & 0<z\le 33+2.5\sin\frac{\pi x}{40}, \\
			7.8, & 33+2.5\sin\frac{\pi x}{40}<z\le 45+0.4x, \\
			7.488, & 45+0.4x<z\le 60+0.4x, \\
			8.268, & 60+0.4x<z\le 100+0.4x, \\
			7.8, & \textnormal{others}.
		\end{array}\right.
	\end{equation*}	
	with unit `km/s'. The model consists of the crust, the mantle and the undulating Moho discontinuity. In the mantle, there is a subduction zone with a thin low velocity layer atop a fast velocity layer \cite{ToYaLiYaHa:16, WuChHuYa:16}, see Figure \ref{fig:pracvel_illu} for illustration. This model is a typical seismogenic zone \cite{ToZhYa:11}. Taking account into the complexity of the velocity structure, it is much more difficult to locate the earthquake hypocenter. The computational time interval $I=[0, 55\,s]$ and the dominant frequency is $f_0=2\,Hz$. Consider $12$ randomly distributed receivers $\boeta_r=(x_r,z_r)$ on the surface $z_r=0\,km$, there horizontal positions are given in Table \ref{tab:pracvel_recepos}.
	
\begin{table*}
    	\caption{The practical velocity model: the horizontal positions of receivers, with unit `km'.} \label{tab:pracvel_recepos}
	\begin{center}\begin{tabular}{ccccccccccccc} \hline
		$r$ & $1$ & $2$ & $3$ & $4$ & $5$ & $6$ & $7$ & $8$ & $9$ & $10$ & $11$ & $12$ \\ \hline
		$x_r$ & $21$ & $33$ & $39$ & $58$ & $68$ & $74$ & $86$ & $98$ & $126$ & $132$ & $158$ & $197$ \\ \hline
    	\end{tabular}\end{center}
\end{table*}

We firstly investigate the noise free situation. Consider the following four cases: (i) An earthquake occurs in the subduction zone, but the initial hypocenter is chosen in the mantle
\begin{align*}
	& \bxi_T=(168.352\,km,142.849\,km),\quad \tau_T=10\,s, \\
	& \bxi=(53.494\,km,47.113\,km),\quad \tau=13.79\,s;
\end{align*}
(ii) The contrary case of (i)
\begin{align*}
	& \bxi_T=(53.494\,km,47.113\,km),\quad \tau_T=10\,s, \\
	& \bxi=(168.352\,km,142.849\,km),\quad \tau=13.79\,s;
\end{align*}
(iii) An earthquake occurs in the crust and close to the Moho discontinuity, but the initial hypocenter is chosen in the subduction zone
\begin{align*}
	& \bxi_T=(163.326\,km,32.877\,km),\quad \tau_T=10\,s, \\
	& \bxi=(26.497\,km,69.235\,km),\quad \tau=15.32\,s;
\end{align*}
(iv) The contrary case of (iii)
\begin{align*}
	& \bxi_T=(26.497\,km,69.235\,km),\quad \tau_T=10\,s, \\
	& \bxi=(163.326\,km,32.877\,km),\quad \tau=15.32\,s;
\end{align*}
The searching domain is $\Omega_s=[0,200\,km]\times[0,200\,km]$, and the searching time interval is $I_s=[\tau-10\,s,\tau+10\,s]$. The mesh sizes for the searching grid are $h_x=0.2\,km$ and $h_z=0.2\,km$. The time step for the searching time interval is $\sigma=0.05\,s$. 

\begin{figure*} 
	\begin{tabular}{ccc}
		\includegraphics[width=0.33\textwidth, height=0.12\textheight]{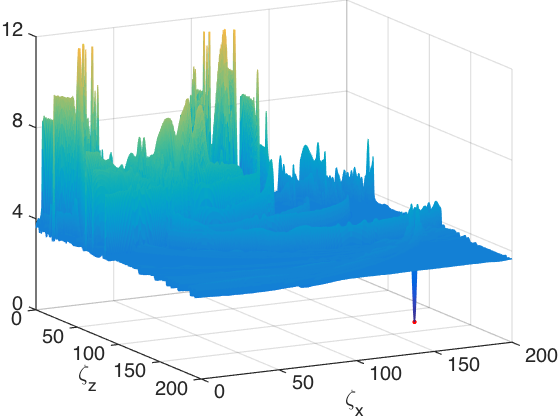} &
		\includegraphics[width=0.33\textwidth, height=0.12\textheight]{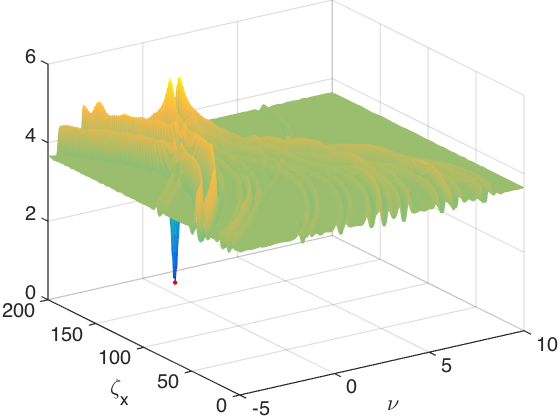} &
		\includegraphics[width=0.33\textwidth, height=0.12\textheight]{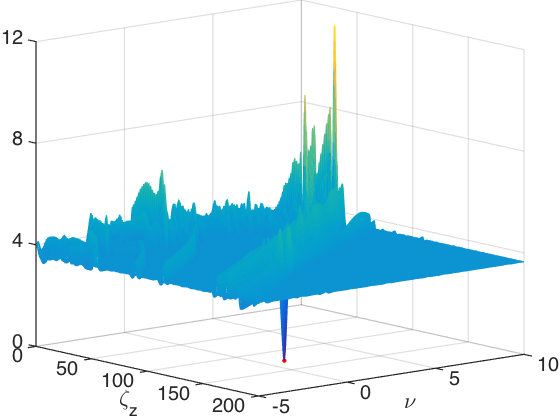} \\
		\includegraphics[width=0.33\textwidth, height=0.12\textheight]{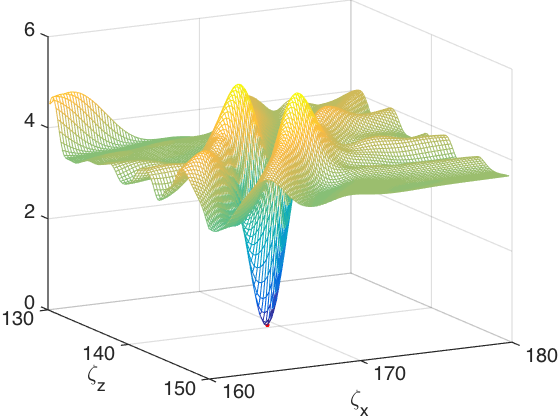} &
		\includegraphics[width=0.33\textwidth, height=0.12\textheight]{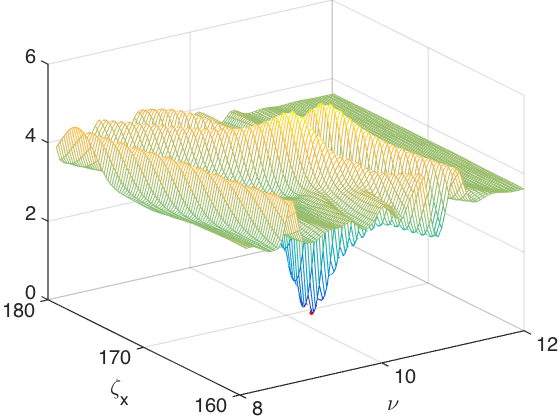} &
		\includegraphics[width=0.33\textwidth, height=0.12\textheight]{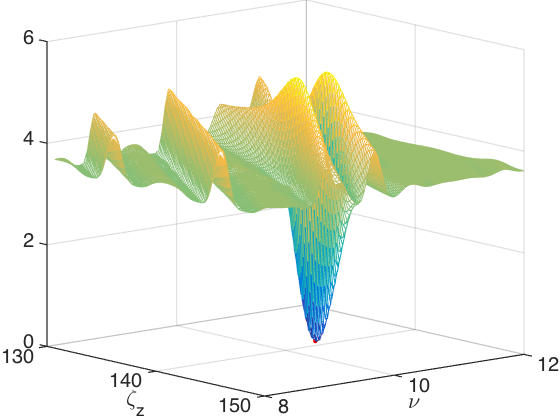}
	\end{tabular}
	\caption{The practical velocity model, case (i). The $\zeta_x-\zeta_z$ (Left, $\nu=10\,s$), $\nu-\zeta_x$ (Middle, $\;\zeta_z=142.849\,km$) and $\nu-\zeta_z$ (Right, $\zeta_x=168.352\,km$) cross section plan of function $\Gamma(\bzeta,\nu)$. Up for large image and Down for zoom-in image.} \label{fig:exam32_4exam_c1auxfun}
\end{figure*}

\begin{figure*} 
	\begin{tabular}{ccc}
		\includegraphics[width=0.33\textwidth, height=0.12\textheight]{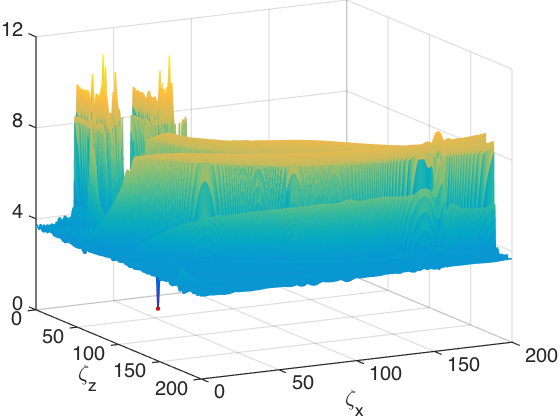} &
		\includegraphics[width=0.33\textwidth, height=0.12\textheight]{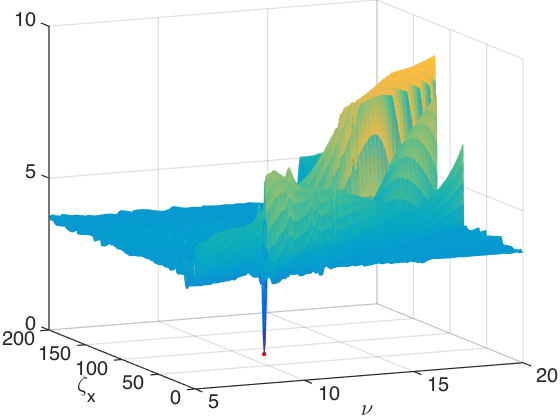} &
		\includegraphics[width=0.33\textwidth, height=0.12\textheight]{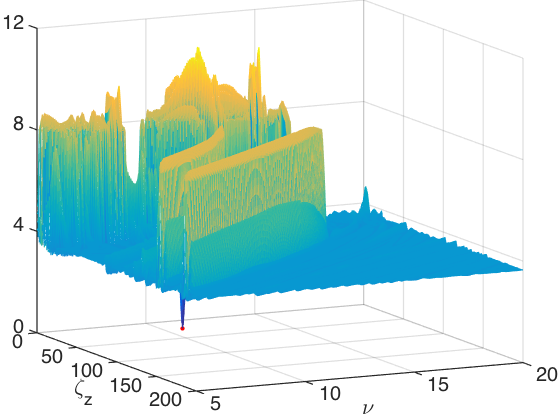} \\
		\includegraphics[width=0.33\textwidth, height=0.12\textheight]{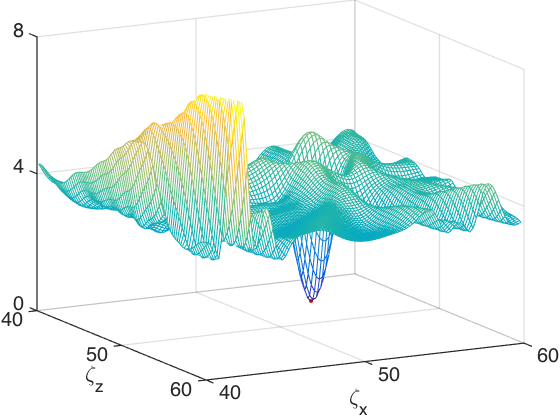} &
		\includegraphics[width=0.33\textwidth, height=0.12\textheight]{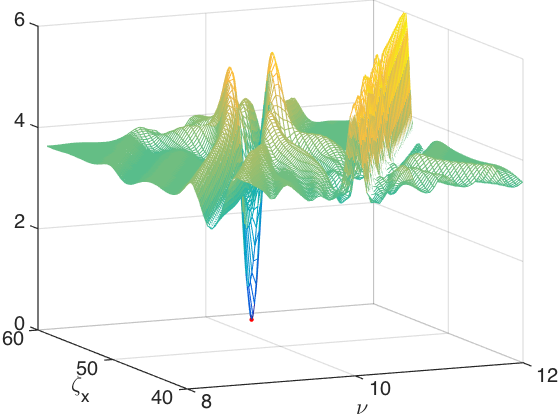} &
		\includegraphics[width=0.33\textwidth, height=0.12\textheight]{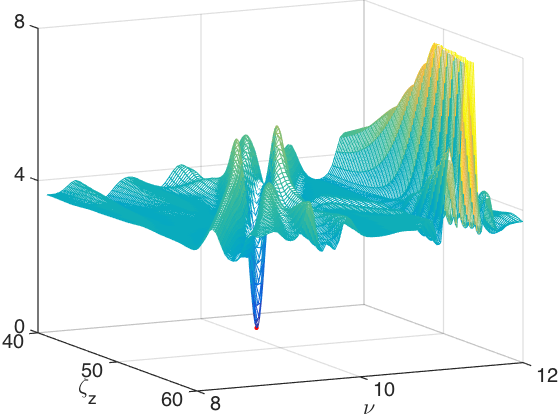}
	\end{tabular}
	\caption{The practical velocity model, case (ii). The $\zeta_x-\zeta_z$ (Left, $\nu=10\,s$), $\nu-\zeta_x$ (Middle, $\;\zeta_z=47.113\,km$) and $\nu-\zeta_z$ (Right, $\zeta_x=53.494\,km$) cross section plan of function $\Gamma(\bzeta,\nu)$. Up for large image and Down for zoom-in image.} \label{fig:exam32_4exam_c2auxfun}
\end{figure*}

\begin{figure*} 
	\begin{tabular}{ccc}
		\includegraphics[width=0.33\textwidth, height=0.12\textheight]{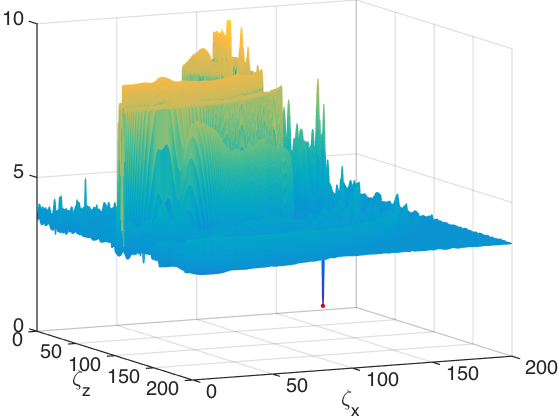} &
		\includegraphics[width=0.33\textwidth, height=0.12\textheight]{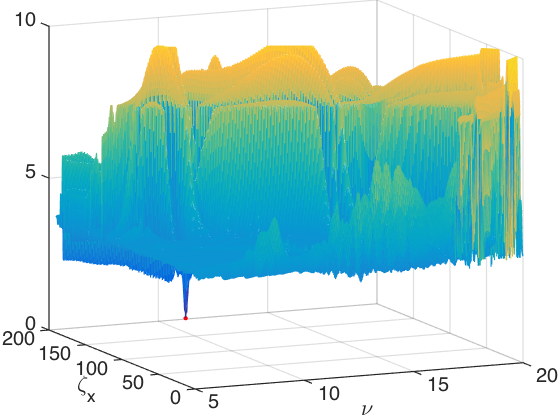} &
		\includegraphics[width=0.33\textwidth, height=0.12\textheight]{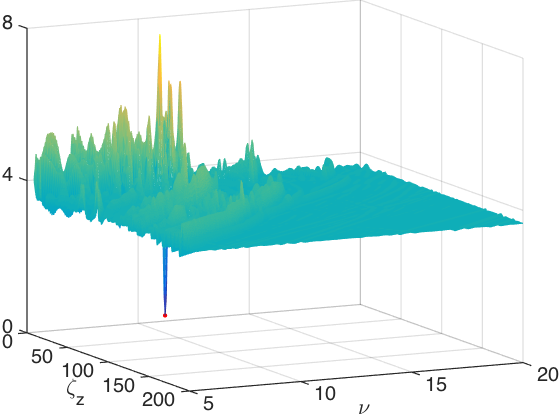} \\
		\includegraphics[width=0.33\textwidth, height=0.12\textheight]{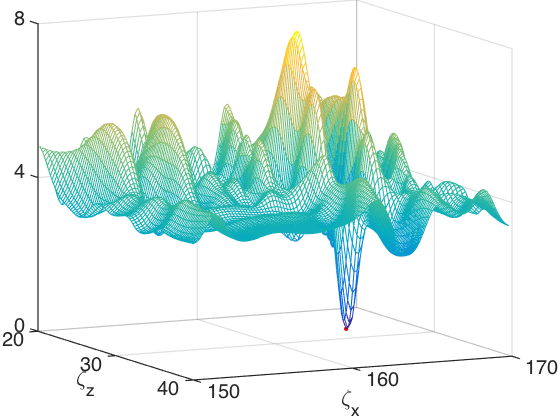} &
		\includegraphics[width=0.33\textwidth, height=0.12\textheight]{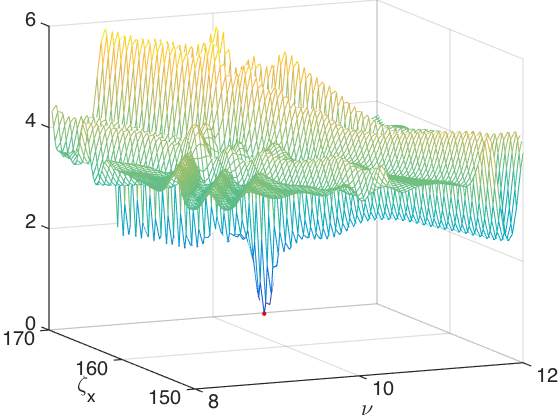} &
		\includegraphics[width=0.33\textwidth, height=0.12\textheight]{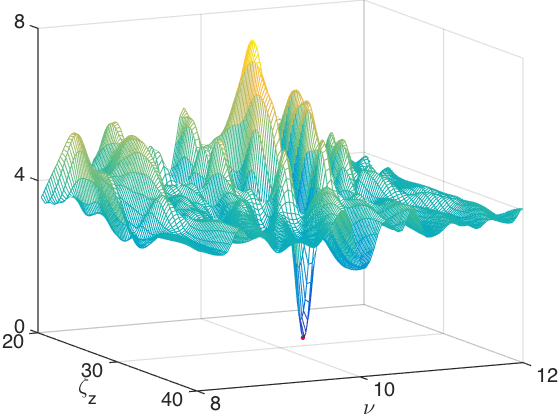}
	\end{tabular}
	\caption{The practical velocity model, case (iii). The $\zeta_x-\zeta_z$ (Left, $\nu=10\,s$), $\nu-\zeta_x$ (Middle, $\;\zeta_z=32.877\,km$) and $\nu-\zeta_z$ (Right, $\zeta_x=163.326\,km$) cross section plan of function $\Gamma(\bzeta,\nu)$. Up for large image and Down for zoom-in image.} \label{fig:exam32_4exam_c3auxfun}
\end{figure*}

\begin{figure*} 
	\begin{tabular}{ccc}
		\includegraphics[width=0.33\textwidth, height=0.12\textheight]{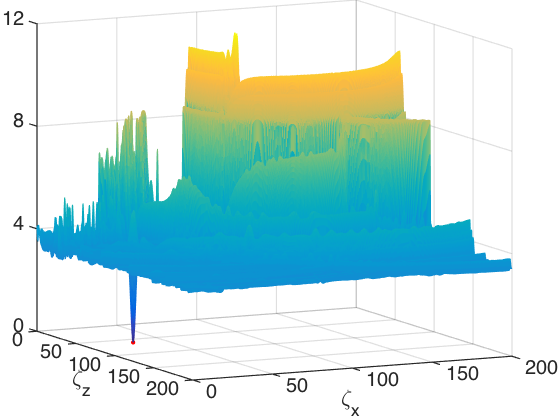} &
		\includegraphics[width=0.33\textwidth, height=0.12\textheight]{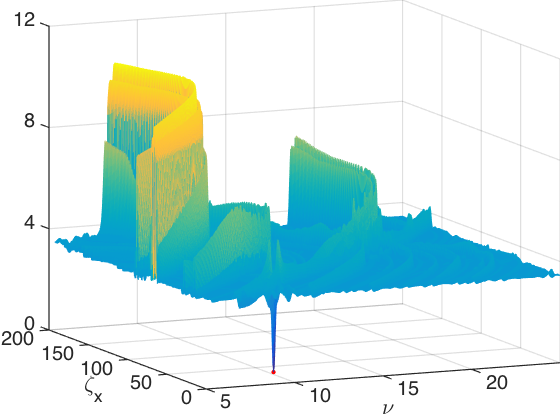} &
		\includegraphics[width=0.33\textwidth, height=0.12\textheight]{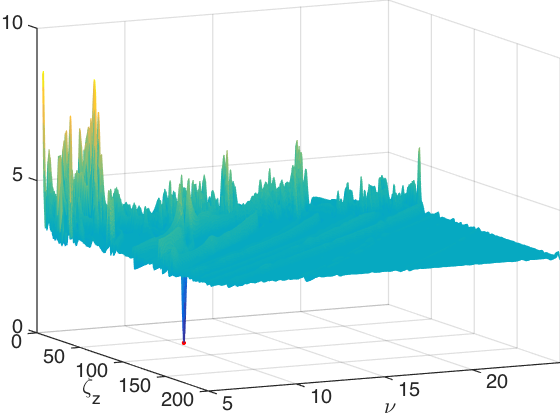} \\
		\includegraphics[width=0.33\textwidth, height=0.12\textheight]{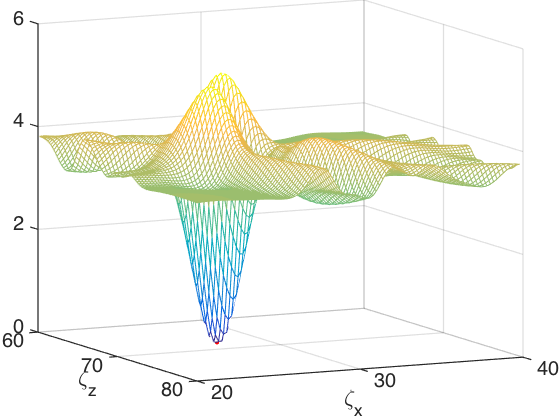} &
		\includegraphics[width=0.33\textwidth, height=0.12\textheight]{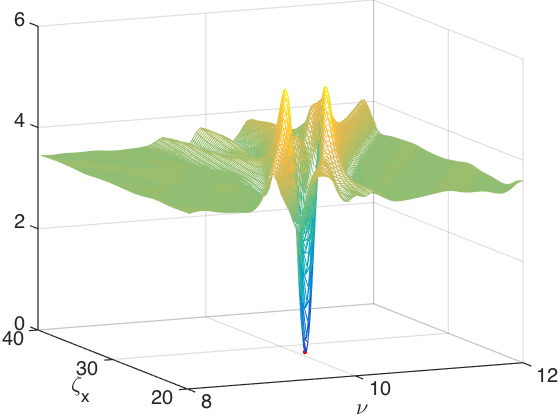} &
		\includegraphics[width=0.33\textwidth, height=0.12\textheight]{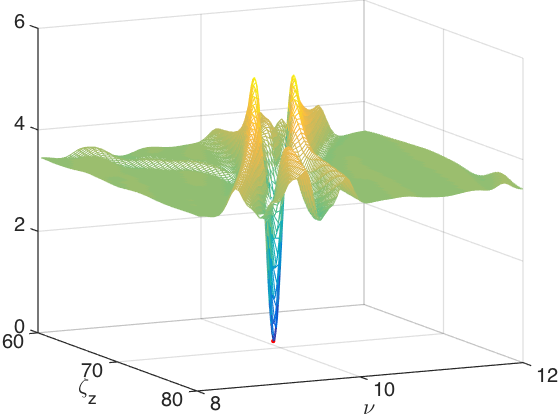}
	\end{tabular}
	\caption{The practical velocity model, case (iv). The $\zeta_x-\zeta_z$ (Left, $\nu=10\,s$), $\nu-\zeta_x$ (Middle, $\;\zeta_z=69.235\,km$) and $\nu-\zeta_z$ (Right, $\zeta_x=26.497\,km$) cross section plan of function $\Gamma(\bzeta,\nu)$. Up for large image and Down for zoom-in image.} \label{fig:exam32_4exam_c4auxfun}
\end{figure*}

\begin{figure*} 
	\begin{tabular}{ccc}
		\includegraphics[width=0.33\textwidth, height=0.12\textheight]{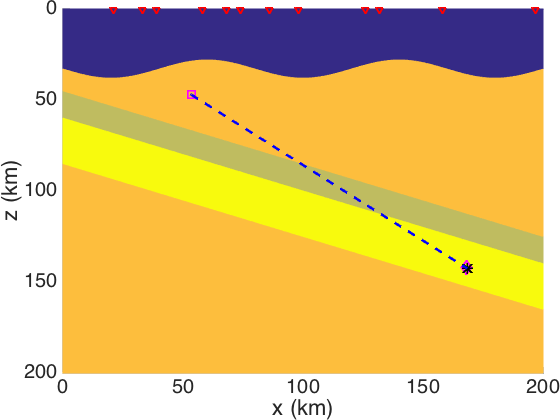} &
		\includegraphics[width=0.33\textwidth, height=0.12\textheight]{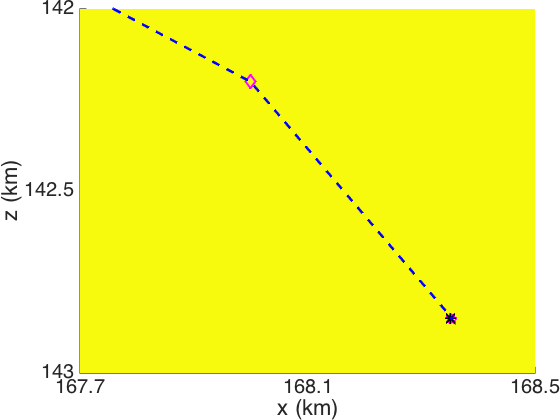} &
		\includegraphics[width=0.33\textwidth, height=0.12\textheight]{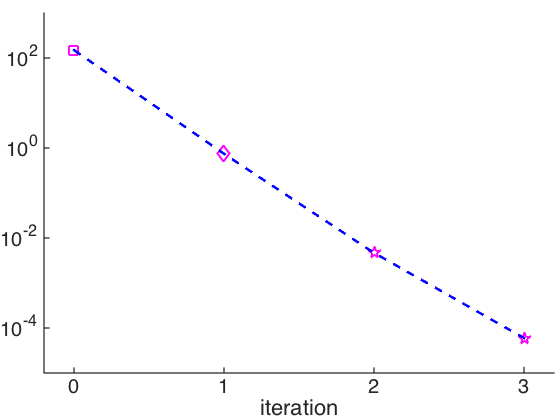} \\
		\includegraphics[width=0.33\textwidth, height=0.12\textheight]{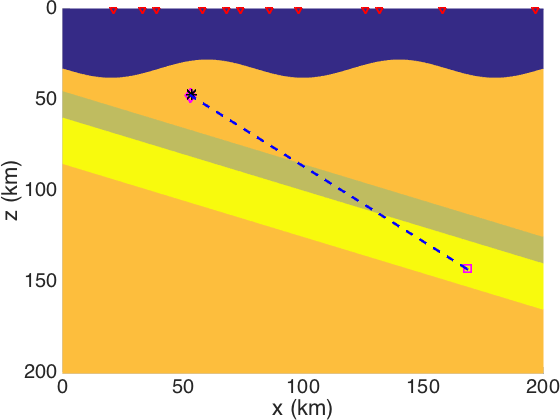} &
		\includegraphics[width=0.33\textwidth, height=0.12\textheight]{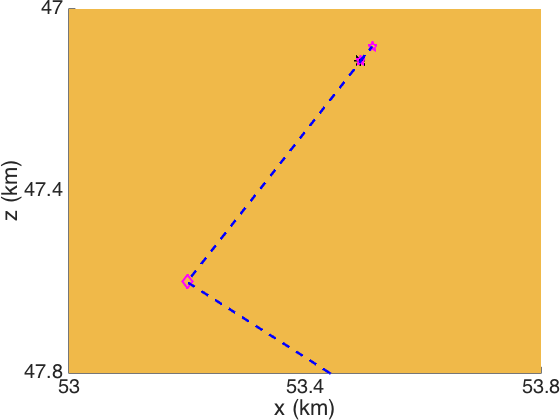} &
		\includegraphics[width=0.33\textwidth, height=0.12\textheight]{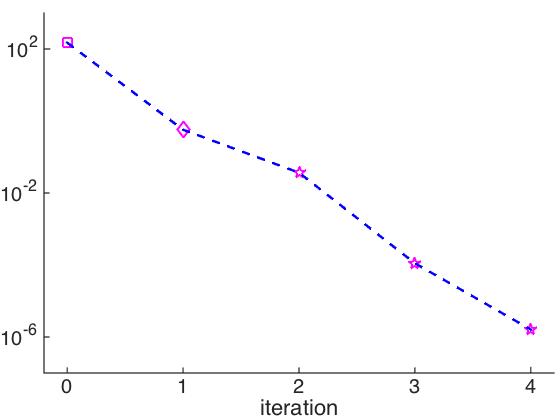} \\
		\includegraphics[width=0.33\textwidth, height=0.12\textheight]{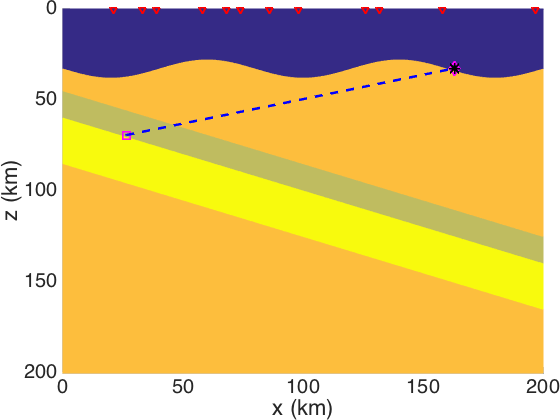} &
		\includegraphics[width=0.33\textwidth, height=0.12\textheight]{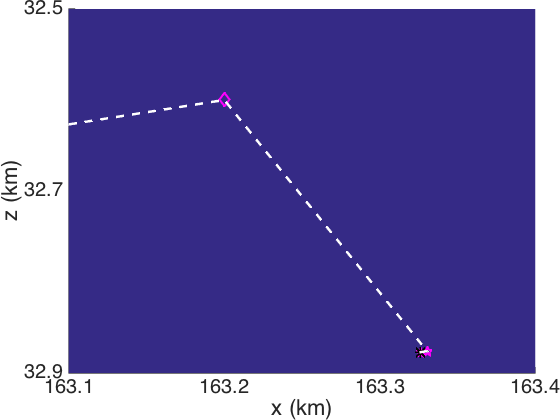} &
		\includegraphics[width=0.33\textwidth, height=0.12\textheight]{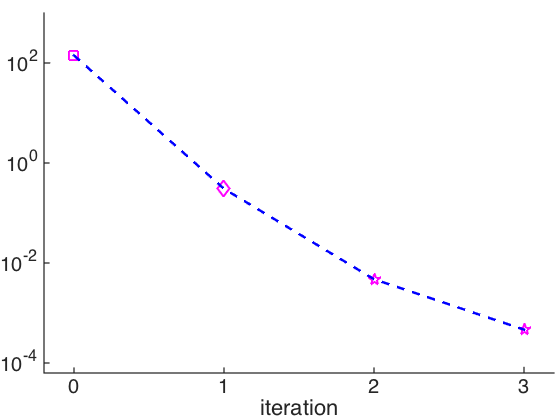} \\
		\includegraphics[width=0.33\textwidth, height=0.12\textheight]{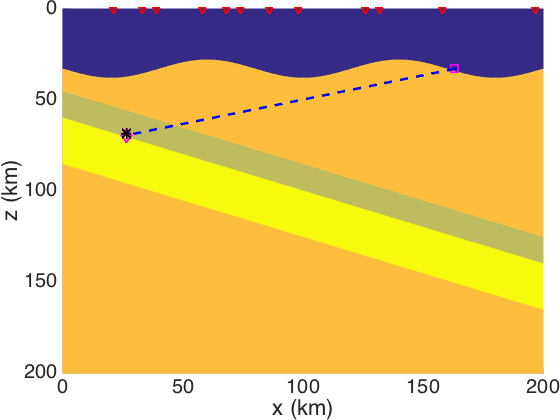} &
		\includegraphics[width=0.33\textwidth, height=0.12\textheight]{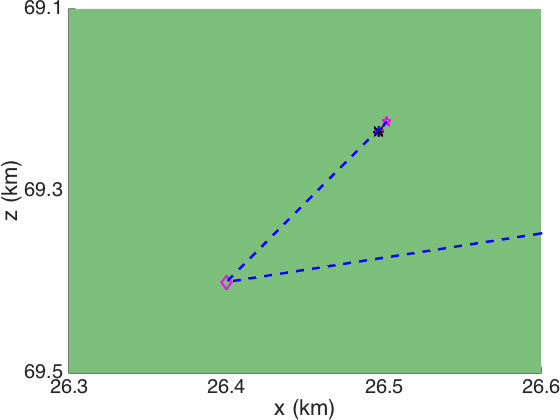} &
		\includegraphics[width=0.33\textwidth, height=0.12\textheight]{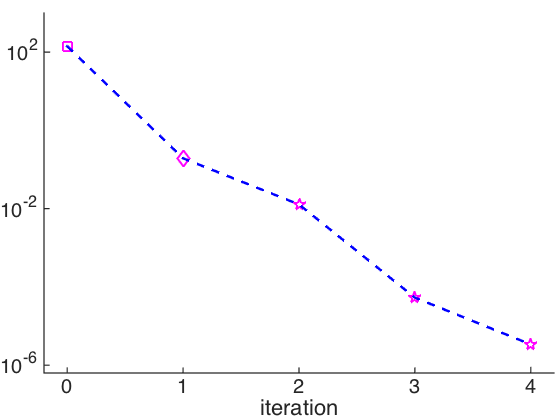}
	\end{tabular}
	\caption{Convergence history of the practical velocity model.  From Up to Down corresponding the cases (i)-(iv). Left (large image) and Middle (zoom-in image): the convergent trajectories; Right: the absolute errors with respect to iteration steps between the real and computed hypocenter of the earthquake. The magenta square is the initial hypocenter, the magenta diamond denotes the hypocenter obtained by AFM (algorithm \ref{alg:afm_class}), the magenta pentagrams indicate the hypocenters in the iterative process of the AFPM (algorithm \ref{alg:afm_preproc}), and the black pentagram is the real hypocenter.} \label{fig:exam32_4exam_trajectory}
\end{figure*}

The functions $\Gamma(\bzeta,\nu)$ are output in Figure \ref{fig:exam32_4exam_c1auxfun}-\ref{fig:exam32_4exam_c4auxfun}. In these figures, it is easy to observe the uniqueness of the global minimum. We also present the convergent history of the auxiliary preprocessor method in Figure \ref{fig:exam32_4exam_trajectory}. We can see that the global minimum of the function $\Gamma(\bzeta,\nu)$ is very close to the optimization solution of \eqref{eqn:optim}. From which, we can draw the same conclusion as in Subsection \ref{subsec:two_layer}.

\begin{figure*} 
	\begin{tabular}{ccc}
		\includegraphics[width=0.33\textwidth, height=0.12\textheight]{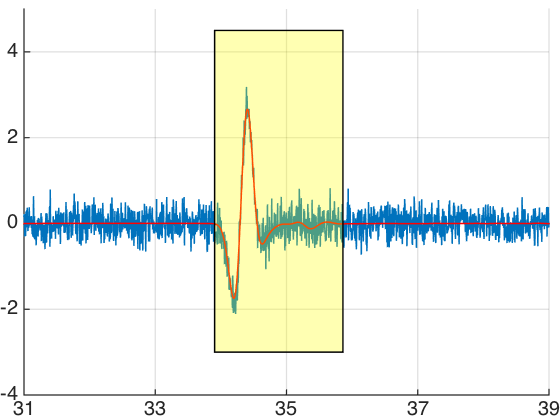} &
		\includegraphics[width=0.33\textwidth, height=0.12\textheight]{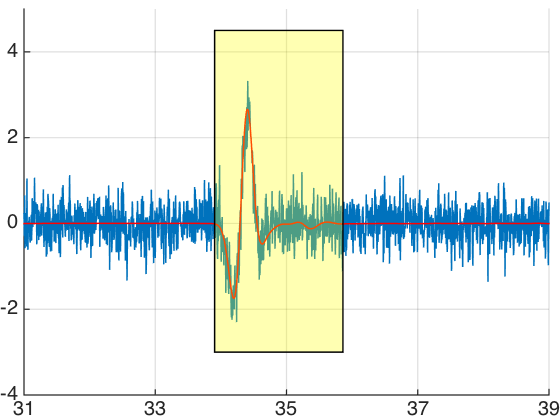} &
		\includegraphics[width=0.33\textwidth, height=0.12\textheight]{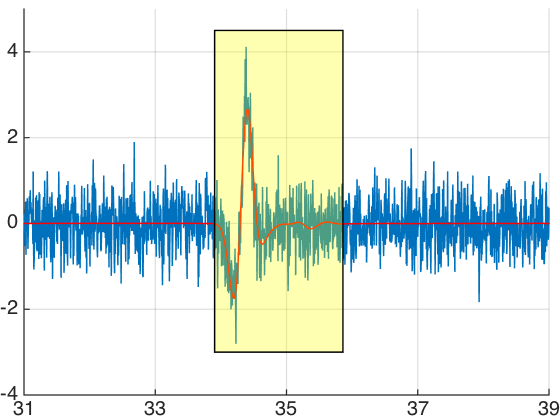} \\
		\includegraphics[width=0.33\textwidth, height=0.12\textheight]{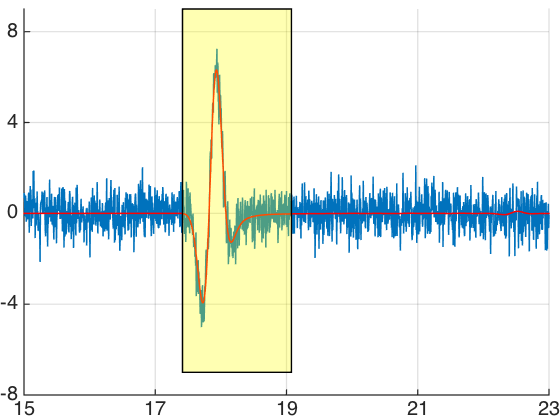} &
		\includegraphics[width=0.33\textwidth, height=0.12\textheight]{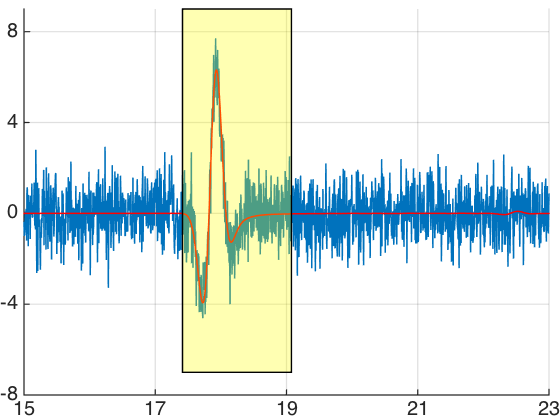} &
		\includegraphics[width=0.33\textwidth, height=0.12\textheight]{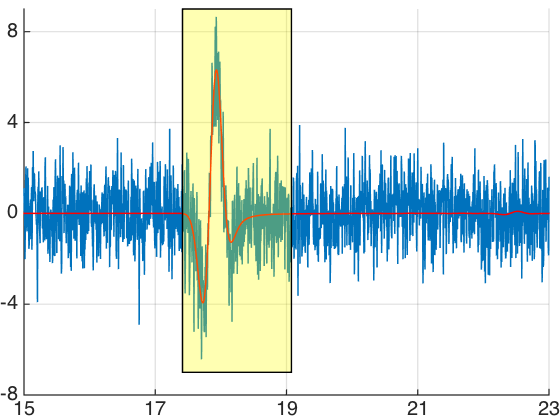} \\
		\includegraphics[width=0.33\textwidth, height=0.12\textheight]{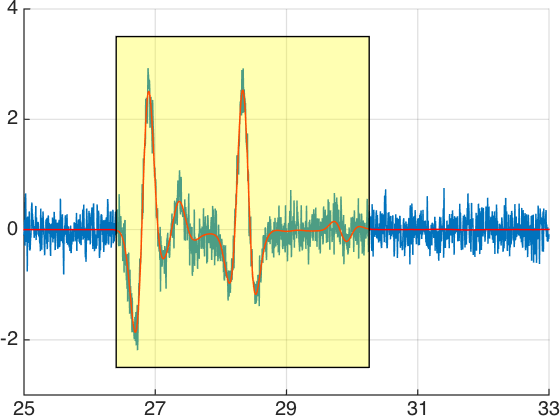} &
		\includegraphics[width=0.33\textwidth, height=0.12\textheight]{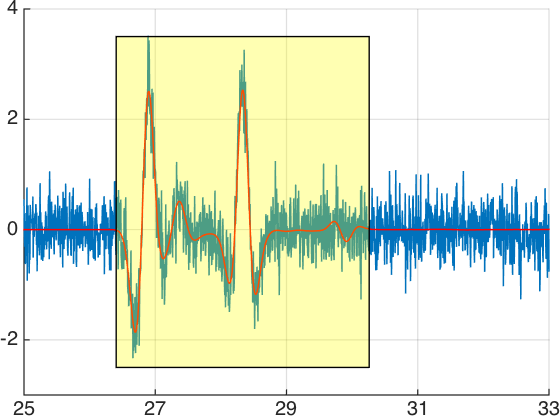} &
		\includegraphics[width=0.33\textwidth, height=0.12\textheight]{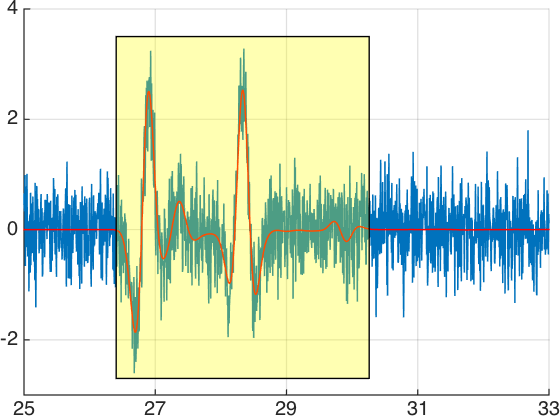} \\
		\includegraphics[width=0.33\textwidth, height=0.12\textheight]{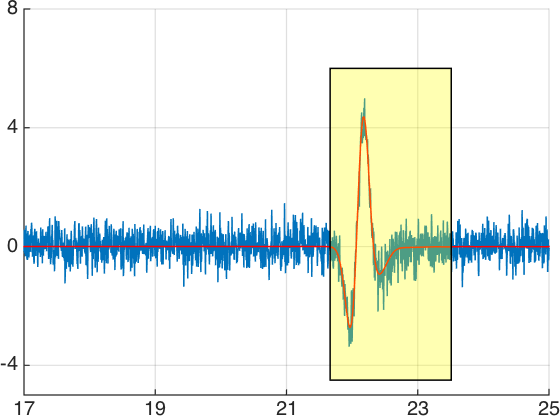} &
		\includegraphics[width=0.33\textwidth, height=0.12\textheight]{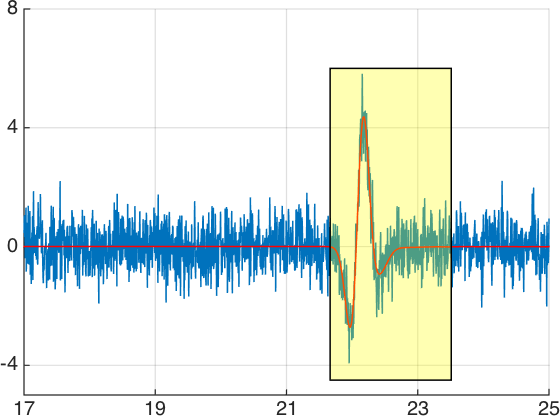} &
		\includegraphics[width=0.33\textwidth, height=0.12\textheight]{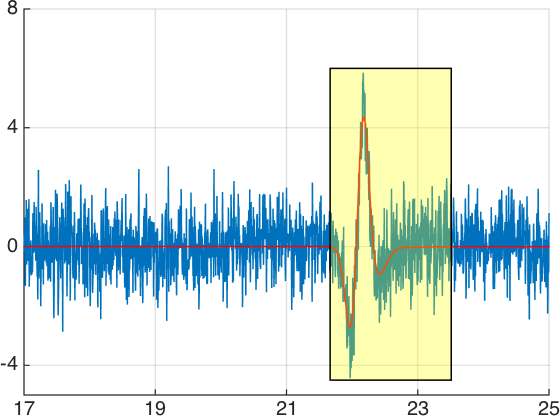}
	\end{tabular}
	\caption{Illustration of signal with noise in the practical velocity model. The signal with noise $d_r(t)$ (blue line) and the noise free signal $u(\boeta_r,t;\bxi_T,\tau_T)$ for receivers $r=5$. The horizontal axis is the time $t$. From Up to Down corresponding the cases (i)-(iv). From left to right, the ratio $R=10\%,\;15\%,\;20\%$ respectively.} \label{fig:exam32_noise}
\end{figure*}

Finally, the noise is taken into consideration. We test the cases (i)-(iv). The noise is added in the same way as in Subsection \ref{subsec:two_layer}. The real earthquake signal with noise $d_r(t)$ and the noise free signal $u(\boeta_r,t;\bxi_T,\tau_T)$ are illustrated in Fig \ref{fig:exam32_noise}. We also select a time window that contains the main part of $u(\boeta_r,t;\bxi_T,\tau_T)$ to reduce the impact of noise. In Figure \ref{fig:exam32_noise} (iii), two waveforms are observed, where the first is the surface wave and the latter is the direct wave. In Figure \ref{fig:exam32_noise} (iv), only one waveform is observed. In fact, since the source location is very close to the discontinuity between the low velocity layer and the fast velocity layer in the subduction zone, the direct wave and the reflected wave arrive almost at the same time. The preceding two cases are very typical. 

In Table \ref{tab:exam32_noise_c1}-\ref{tab:exam32_noise_c4}, the location results $(\bxi_*,\tau_*)$ computed via the auxiliary function method (Algorithm \ref{alg:afm_class}) and their errors with respected to different ratio $R$ are presented. From which, we can see that the method can obtain satisfactory location results for $R=10\%$ and $15\%$。For $R=20\%$, the error results are obtained for case (i) and (iii). This is because the simulation domain is very large and there are many discontinuities. Nevertheless, the auxiliary function method is still much better than the iterative methods, which is only valid for $R\le 5\%$ according to our numerical tests. Taking into account that the computational cost of the auxiliary function method is almost the same to the single iteration step of the iterative method. The computation efficiency of our method is also obvious.

\begin{table*}
    	\caption{The practical velocity model, case (i). The location results $(\bxi_*,\tau_*)$ and the errors.} \label{tab:exam32_noise_c1}
	\begin{center}\begin{tabular}{c|cccc} \hline
		 $R$ & $\bxi_*$\;(km) & $\tau_*$\;(s) & $\norm{\bxi_T-\bxi_*}_2$\;(km) & $\abs{\bxi_{*z}-\bxi_{Tz}}/\bxi_{Tz}$ \\ \hline
		 $10\%$ & $(168.2,142.6)$ & $10.04$ & $0.292$ & $0.17\%$ \\
		 $15\%$ & $(168.2,142.2)$ & $10.09$ & $0.667$ & $0.45\%$ \\ \hline
		 $20\%$ & \multicolumn{4}{c}{fail} \\ \hline
    	\end{tabular}\end{center}
\end{table*}

\begin{table*}
    	\caption{The practical velocity model, case (ii). The location results $(\bxi_*,\tau_*)$ and the errors.} \label{tab:exam32_noise_c2}
	\begin{center}\begin{tabular}{c|cccc} \hline
		 $R$ & $\bxi_*$\;(km) & $\tau_*$\;(s) & $\norm{\bxi_T-\bxi_*}_2$\;(km) & $\abs{\bxi_{*z}-\bxi_{Tz}}/\bxi_{Tz}$ \\ \hline
		 $10\%$ & $(53.4,47.2)$ & $9.99$ & $0.128$ & $0.18\%$ \\
		 $15\%$ & $(53.4,47.2)$ & $9.99$ & $0.128$ & $0.18\%$ \\
		 $20\%$ & $(53.4,47.2)$ & $9.99$ & $0.128$ & $0.18\%$ \\ \hline
    	\end{tabular}\end{center}
\end{table*}

\begin{table*}
    	\caption{The practical velocity model, case (iii). The location results $(\bxi_*,\tau_*)$ and the errors.} \label{tab:exam32_noise_c3}
	\begin{center}\begin{tabular}{c|cccc} \hline
		 $R$ & $\bxi_*$\;(km) & $\tau_*$\;(s) & $\norm{\bxi_T-\bxi_*}_2$\;(km) & $\abs{\bxi_{*z}-\bxi_{Tz}}/\bxi_{Tz}$ \\ \hline
		 $10\%$ & $(163.2,32.6)$ & $10.02$ & $0.304$ & $0.84\%$ \\
		 $15\%$ & $(163.2,32.6)$ & $10.02$ & $0.304$ & $0.84\%$\\ \hline
		 $20\%$ & \multicolumn{4}{c}{fail} \\ \hline
    	\end{tabular}\end{center}
\end{table*}

\begin{table*}
    	\caption{The practical velocity model, case (iv). The location results $(\bxi_*,\tau_*)$ and the errors.} \label{tab:exam32_noise_c4}
	\begin{center}\begin{tabular}{c|cccc} \hline
		 $R$ & $\bxi_*$\;(km) & $\tau_*$\;(s) & $\norm{\bxi_T-\bxi_*}_2$\;(km) & $\abs{\bxi_{*z}-\bxi_{Tz}}/\bxi_{Tz}$ \\ \hline
		 $10\%$ & $(26.4,69.4)$ & $9.97$ & $0.191$ & $0.24\%$ \\
		 $15\%$ & $(26.4,69.4)$ & $9.97$ & $0.191$ & $0.24\%$ \\
		 $20\%$ & $(26.4,69.4)$ & $9.97$ & $0.191$ & $0.24\%$ \\ \hline
    	\end{tabular}\end{center}
\end{table*}

\section{Conclusion and Discussion} \label{sec:con}
The first conclusion to be drawn from the numerical evidence presented earlier is that the auxiliary function preprocessor method (Algorithm \ref{alg:afm_preproc}) can determine the earthquake hypocenter and the origin time very efficient and accurate when the seismic signals are noise-free. Secondly, the auxiliary function method (Algorithm \ref{alg:afm_class}) can locate the earthquake hypocenter and the origin time with reasonable accuracy in the situation of noise. The above advantages are based on the fact that the real hypocenter and origin time is the root of the new constructed auxiliary functions. And the total computational cost of constructing these functions is comparable to the single iteration step of the iterative method.

It should be noted that there are still many issues need to be further investigated: (a) Currently, the uniqueness of the solution can not be proved, but we have observed the uniqueness of the solution numerically. If is of course very interesting to present an intuitive study.  (b) We are currently working on the 2-D problem and the accurate velocity model. For 3-D problem and the inaccurate velocity model, the situations may be more complicated. These all require much more effort.


\section*{Acknowledgement}
This work was supported by the National Nature Science Foundation of China (Grant Nos 41230210 and 41390452). Hao Wu was also partially supported by the National Nature Science Foundation of China (Grant Nos 11101236 and 91330203) and SRF for ROCS, SEM. The authors are grateful to Prof. Shi Jin for his helpful suggestions and discussions that greatly improve the presentation. Hao Wu would like to thank Prof. Jingwei Hu and Dr. Meixia Wang for their valuable comments.


\begin{thebibliography}{}
	\bibitem{AkRi:80}
		K. Aki and P.G. Richards, \textit{Quantitative Seismology: Theory and Methods volume II}, W.H. Freeman \& Co (Sd), 1980.
		
	\bibitem{Da:86}
		M.A. Dablain, The application of high-order differencing to the scalar wave equation, \textit{Geophysics}, \textbf{51}(1), 54-66, 1986.		

	\bibitem{EnRu:03}
    		B. Engquist and O. Runborg, Computational high frequency wave propagation, \textit{Acta Numer.}, \textbf{12}, 181-266, 2003.	

	\bibitem{Ge:03}
    		M.C. Ge, Analysis of source location algorithms Part I: Overview and non-iterative methods, \textit{J. Acoust. Emiss.}, \textbf{21}, 14-28. 2003.

	\bibitem{Ge:03b}
    		M.C. Ge, Analysis of source location algorithms Part II: Iterative methods, \textit{J. Acoust. Emiss.}, \textbf{21}, 29-51, 2003.
		
	\bibitem{Ge:12}
    		L. Geiger, Probability method for the determination of earthquake epicenters from the arrival time only, \textit{Bull. St. Louis Univ.}, \textbf{8}, 60-71, 1912.
		
	\bibitem{HuYaToBaQi:16}
    		X.Y. Huang, D.H. Yang, P. Tong, J. Badal and Q.Y. Liu, Wave equation-based reflection tomography of the 1992 Landers earthquake area, \textit{Geophys. Res. Lett.}, \textbf{43}, 1884-1892, 2016.
			
	\bibitem{JiWuYa:08}
    		S. Jin, H. Wu and X. Yang, Gaussian Beam Methods for the Schr\"odinger Equation in the Semi-classical Regime: Lagrangian and Eulerian Formulations, \textit{Commun. Math. Sci.}, \textbf{6}(4), 995-1020, 2008.
				
	\bibitem{KiLiTr:11}
		Y.H. Kim, Q.Y. Liu and J. Tromp, Adjoint centroid-moment tensor inversions, \textit{Geophys. J. Int.}, \textbf{186}, 264-278, 2011.

	\bibitem{KoTr:03}
		D. Komatitsch and J. Tromp, A perfectly matched layer absorbing boundary condition for the second-order seismic wave equation, \textit{Geophys. J. Int.}, \textbf{154}, 146-153, 2003.
		
	\bibitem{KoTr:03}
		D. Komatitsch and J. Tromp, A perfectly matched layer absorbing boundary condition for the second-order seismic wave equation, \textit{Geophys. J. Int.}, \textbf{154}, 146-153, 2003.
		
	\bibitem{LeSt:81}
		W.H.K. Lee and S.W. Stewart, \textit{Principles and Applications of Microearthquake Networks}, Academic Press, 1981.
		
	\bibitem{LiGu:12}
		Q.Y. Liu and Y.J. Gu, Seismic imagine: From classical to adjoint tomography, \textit{Tectonophysics}, \textbf{566-567}, 31-66, 2012.

	\bibitem{LiPoKoTr:04}
		Q.Y. Liu, J. Polet, D. Komatitsch and J. Tromp, Spectral-Element Moment Tensor Inversion for Earthquakes in Southern California, \textit{Bull. seism. Soc. Am.}, \textbf{94}(5), 1748-1761, 2004.
		
	\bibitem{Ma:15}
		R. Madariaga, Seismic Source Theory, in \textit{Treatise on Geophysics (Second Edition)}, pp. 51-71, ed. Gerald, S., Elsevier B.V., 2015.

	\bibitem{PrGe:88}
		A.F. Prugger and D.J. Gendzwill, Microearthquake location: A nonlinear approach that makes use of a simplex stepping procedure, \textit{Bull. seism. Soc. Am.}, \textbf{78}, 799-815, 1988.
		
	\bibitem{RaPoFi:10}
    		N. Rawlinson, S. Pozgay and S. Fishwick, Seismic tomography: A window into deep Earth, \textit{Phys. Earth Planet. Inter.}, \textbf{178}, 101-135, 2010.
				
	\bibitem{SaLoZo:08}
		C. Satriano, A. Lomax and A. Zollo, Real-Time Evolutionary Earthquake Location for Seismic Early Warning, \textit{Bull. seism. Soc. Am.}, \textbf{98}(3), 1482-1494, 2008.
		
	\bibitem{Th:85}
		C.H. Thurber, Nonlinear earthquake location: Theory and examples, \textit{Bull. seism. Soc. Am.}, \textbf{75}(3), 779-790, 1985.

	\bibitem{Th:14}
		C.H. Thurber, Earthquake, location techniques, in \textit{Encyclopedia of Earth Sciences Series}, pp. 201-207, ed. Gupta, H.K., Springer, 2014.
		
	\bibitem{ToZhYa:11}
		P. Tong, D.P. Zhao and D.H. Yang, Tomography of the 1995 Kobe earthquake area: comparison of finite-frequency and ray approaches, \textit{Geophys. J. Int.}, \textbf{187}, 278-302, 2011.		
	\bibitem{ToZhYaYaChLi:14}
		P. Tong, D. Zhao, D.H. Yang, X. Yang, J. Chen and Q. Liu, Wave-equation-based travel-time seismic tomography - Part 1: Method, \textit{Solid Earth}, \textbf{5}, 1151-1168, 2014.
		
	\bibitem{ToYaLiYaHa:16}
		P. Tong, D.H. Yang, Q.Y. Liu, X. Yang and J. Harris, Acoustic wave-equation-based earthquake location, \textit{Geophys. J. Int.}, \textbf{205}(1), 464-478, 2016.
				
	\bibitem{WaEl:00}
		F. Waldhauser and W.L. Ellsworth, A double-difference earthquake location algorithm: Method and application to the northern Hayward Fault, California, \textit{Bull. seism. Soc. Am.}, \textbf{90}(6), 1353-1368, 2000.
		
	\bibitem{We:08}
    		X. Wen, High Order Numerical Quadratures to One Dimensional Delta Function Integrals, \textit{SIAM J. Sci. Comput.}, \textbf{30}(4), 1825-1846, 2008.

	\bibitem{WuChHuYa:16}
    		H. Wu, J. Chen, X.Y. Huang and D.H. Yang, A new earthquake location method based on the waveform inversion, \textit{Commun. Comput. Phys.}, accepted, 2017.
				
	\bibitem{WuYa:13}
    		H. Wu and X. Yang, Eulerian Gaussian beam method for high frequency wave propagation in the reduced momentum space, \textit{Wave Motion}, \textbf{50}(6), 1036-1049, 2013.

	\bibitem{YaLuWuPe:04}
		D.H. Yang, M. Lu, R.S. Wu and J.M. Peng, An Optimal Nearly Analytic Discrete Method for 2D Acoustic and Elastic Wave Equations, \textit{Bull. seism. Soc. Am.}, \textbf{94}(5), 1982-1991, 2004.

\end{thebibliography}
\end{document}